\documentclass[preprint, review, 12pt]{elsarticle}

\usepackage{amssymb}
\usepackage{amsthm}
\usepackage{fullpage}
\usepackage{graphicx}
\usepackage{multicol}
\usepackage{marvosym}
\usepackage{tikz}
\usetikzlibrary{matrix,calc}
\usetikzlibrary{patterns}
\usepackage{tensor}
\usepackage{amsmath}
\usepackage{tkz-euclide}
\usepackage{subfigure}
\usetkzobj{all}
\usepackage[colorlinks = true, linkcolor = black, citecolor = black, final]{hyperref}

\usepackage{amsthm}
 
\theoremstyle{definition}
\newtheorem{definition}{Definition}[section]

\journal{Journal Name}

\begin{document}

\begin{frontmatter}

\title{3rd-order Spectral Representation Method: Part II -- Ergodic Multi-variate random processes with fast Fourier transform}

\author{Lohit Vandanapu}
\author{Michael D. Shields}

\address{Department of Civil Engineering, Johns Hopkins University}

\begin{abstract}
The second in a two-part series, this paper extends the 3rd-order Spectral Representation Method for simulation of ergodic multi-variate stochastic processes according to a prescribed cross power spectral density and cross bispectral density. The $2^{nd}$ and $3^{rd}$ order ensemble properties of the simulated stochastic vector processes are shown to satisfy the target cross correlation properties in expectation. A multi-indexed frequency discretization is introduced to ensure ergodicity of the sample functions. This is first shown for uni-variate processes and then the simulation formula for multi-variate processes is provided. Ensemble properties and ergodicity of the sample functions are proven. Additionally, it is shown that the simulations can be implemented efficiently with the Fast Fourier Transform, which greatly reduces computational effort.  An example involving the simulation of turbulent wind velocity fluctuations is presented to further highlight the features and applications of the algorithm.
\end{abstract}


\end{frontmatter}

\section{Introduction}
\label{S:1}

Numerous civil and mechanical systems involve uncertainties that can be characterized by stochastic processes or fields. Even multiples decades after its inception, monte Carlo Simulation remains the benchmark method for solving a wide-range of non-linear stochastic mechanics problems. In these cases, the monte Carlo simulation framework requires generation of sample functions of the stochastic processes and fields. A brief review of methods used for the simulation of stochastic processes/fields is presented in Part I of this article \cite{VandanapuI}.

In this work, we are specifically interested in the simulation of multi-variate stochastic processes (stochastic vector processes), which are composed of multiple correlated random processes occurring at different discrete spatial locations. Several methods, dating back nearly 50 years, have been proposed for the simulation of stochastic vector processes including methods for stationary, non-stationary, Gaussian and non-Gaussian processes. much of this began with the seminal work of Shinozuka who proposed the Spectral Representation method (SRM) in the early 1970s \cite{shinozuka1971simulation}. Later, in the late 1980s, Mignolet and Spanos \cite{mignolet1987, Spanos1987}, in a 2-part paper, introduced the recursive simulation of stationary multivariate stochastic processes based on autoregressive moving averages methods. This was followed by numerous works in the 1990s, when much of the theory for the SRM was developed. \citet{Li1991} developed a framework for the simulation of non-stationary multi-variate processes with the use of a stochastic decomposition technique and later developed a hybrid discrete Fourier Transform and digital filtering approach \cite{Li1993}. With regard to the SRM, Deodatis developed the theoretical framework for simulation of ergodic, Gaussian stochastic vector processes in 1996 \cite{Deodatis1996}, with subsequent extensions to non-stationary \cite{deodatis1996non} and non-Gaussian processes \cite{popescu1998simulation}.  more recently, the iterative translation approximation method (ITAM) has been proposed for the efficient simulation of non-Gaussian stochastic vector translation processes by Shields and Deodatis \cite{shields2013simple}. Very recently, \citet{Liu2018}, proposed a novel method based on the combination of SRM with a proper orthogonal decomposition for dimension reduction.

The purpose of this Part II article is to extend the 3rd-order Spectral Representation method developed in \cite{Shields2017, VandanapuI} to simulate non-Gaussian multi-variate random processes with specified cross-spectral density and cross-bispectral density. Higher-order (non-Gaussian) multi-variate random processes occur often in nature when non-Gaussian random processes at different spatial locations are related through a specified correlation. For example, seismic acceleration records or wind velocity/pressure time histories at nearby locations may be non-Gaussian in nature and strongly correlated. As such, it is important to model these spatially correlated random processes for use with physics-based models of, for example, structures subjected to wind, seismic, or sea wave excitations at multiple locations. 

Another important contribution of the work is the introduction of a modified third-order SRM simulation formula for the generation of ergodic stochastic processes. The existing formulation presented in \cite{Shields2017, VandanapuI} produces realizations that are non-ergodic and a doubled frequency indexing is presented herein to ensure ergodicity of the sample functions. Finally, we show that the simulation formula can be implemented with the powerful Fast Fourier Transform (FFT), which leads to tremendous gains in computational time.

As a notational note, this Part II paper differs from Part I in the following ways. Since Part I focused on multi-dimensional processes whose randomness is indexed on a spatial variable, random fields were indexed on the spatial variable $x$ with spatial separation denoted by $\xi$ and wave number $\kappa$. This Part II paper, on the other hand, deals with temporal randomness and therefore the random processes are indexed on the time variable $t$ with time lag $\tau$ and circular frequency denoted by $\omega$.

\section{Properties of Stochastic Vector Processes}
\label{S:2}

Before developing the simulation methodologies, we present several important definitions and properties of stochastic vector processes. First, the concepts of stationarity (up to various orders) are presented. Next, the ensemble properties of multi-variate processes are presented. Finally, the spectral properties are discussed.


Let us first begin with a basic definition of a stochastic vector process:
\theoremstyle{definition}
\begin{definition}{}
A complex-valued random vector $\boldsymbol{f}$ is considered a \textbf{stochastic vector process} if its individual individual components $\{ f_{i}\}$ are all indexed by the same continuous time parameter $t$ as:
\begin{equation}
\begin{aligned}
	&\boldsymbol{f}(t) = [f_{1}(t), f_{2}(t), \dots, f_{m}(t)].
\end{aligned}
\end{equation}
\end{definition}

\subsection{Stationary Stochastic Vector Processes}

Stationarity refers to the probabilistic invariance of the random process/field under a shift in the indexing parameter (time or space). Similar to Part I, we present 3 notions of stationarity for stochastic vector processes, although they are presented somewhat more succinctly here.

\subsubsection{Strictly or Strongly Stationary Random Processes}
\label{S:strict_stationarity}

A stochastic vector process $\boldsymbol{f}(t)$ is said to be strictly stationary, or strongly stationary, if the complete joint probability measure is invariant to a shift in index. For instance, suppose that $\boldsymbol{f}(t)$ has $m\times M$-dimensional joint cumulative distribution function for any finite number of stochastic vector processes $m$ having $m$ discrete time increments, given by $F(f_1(t_{1}), f_1(t_{2}),\dots, f_1(t_{m}), \dots, f_m(t_{1}), f_m(t_{2}),\dots, f_m(t_{m}))$, the vector process is said to be strongly stationary if
\begin{equation}
	\begin{aligned}
	& F(f_1(t_{1}), f_1(t_{2}),\dots, f_1(t_{M}), \dots, f_m(t_{1}), f_m(t_{2}),\dots, f_m(t_{M})) =\\
	& F(f_1(t_{1}+\tau), f_1(t_{2}+\tau),\dots, f_1(t_{M}+\tau), \dots, f_m(t_{1}+\tau), f_m(t_{2}+\tau),\dots, f_m(t_{M}+\tau)),\quad \forall \tau, 
	\end{aligned}
    \label{eqn:strong_stationary}
\end{equation}
It follows directly that all lower-dimensional distributions are similarly invariant to a shift in index, $\tau$, and that any characteristic of the joint distribution (i.e.\ moments, cumulants, etc.) are independent of $\tau$.



\subsubsection{$k^{th}$-order Stationary Random Processes}

Following the definitions presented in \ref{S:strict_stationarity}, a stochastic vector process is said to be $k^{th}$ order stationary if the process is probabilistically invariant up to order $k$. That is, the following condition is satisfied:
\begin{equation}
	\begin{aligned}
	& F^{(k)}(f_1(t_{1}), f_1(t_{2}),\dots, f_1(t_{k}), \dots, f_m(t_{1}), f_m(t_{2}),\dots, f_m(t_{k})) =\\
	& F^{(k)}(f_1(t_{1}+\tau), f_1(t_{2}+\tau),\dots, f_1(t_{k}+\tau), \dots, f_m(t_{1}+\tau), f_m(t_{2}+\tau),\dots, f_m(t_{k}+\tau)),\quad \forall \tau, 
	\end{aligned}
    \label{eqn:strong_stationary}
\end{equation}
where $F^{(k)}(\cdot)$ is the $k^{th}$ order joint cumulative distribution function. It is again apparent that all measures of order $< k$ are similarly invariant to a shift in the index $\tau$. It is important to understand the concepts of $3^{rd}$-order stationarity in stochastic vector processes as the proposed methodology simulates vector processes that are $3^{rd}$-order stationary.


\subsubsection{Weak or Wide-Sense Stationary Random Processes}

A random process is considered to be weakly, or wide-sense stationary if the joint probability distribution up to $2^{nd}$-order is invariant to a shift in index. In other words, a weakly stationary random process is a $k^{th}$-order random process with $k=2$. 

\subsection{Ensemble Properties of Stochastic Vector Processes}

Consider a one-dimensional, m-variate (1D-mV) third-order stationary stochastic vector process $\boldsymbol{f}(t)$ with components $[f_{1}(t), f_{2}(t), ....... f_{m}(t)]$ having zero mean for each component,
\begin{equation}
    \mathbb{E}[f_{j}(t)] = 0 \ \text{for} \ j=1, 2, \dots, m
\end{equation}
The second-order correlation function of this stochastic vector process is given by
\begin{equation}
    \mathbb{E}[f_{i}(t)f_{j}(t + \tau)] = R_{ij}(\tau) \ \text{for} \ i, j=1, 2, \dots, m
\end{equation}
which is generally provided through the cross-correlation matrix given by
\begin{equation}
R(\tau) = 
\begin{bmatrix}
    R_{11}(\tau) & R_{12}(\tau) & R_{13}(\tau) & \dots  & R_{1m}(\tau) \\
    R_{21}(\tau) & R_{22}(\tau) & R_{23}(\tau) & \dots  & R_{2m}(\tau) \\
    R_{31}(\tau) & R_{32}(\tau) & R_{33}(\tau) & \dots  & R_{3m}(\tau) \\
    \vdots        & \vdots        & \vdots        & \ddots & \vdots        \\
    R_{m1}(\tau) & R_{m2}(\tau) & R_{m3}(\tau) & \dots  & R_{mm}(\tau) \\
\end{bmatrix}
\end{equation}

Similarly, the third-order correlation function of a stochastic vector process can be expressed as follows:
\begin{equation}
    \mathbb{E}[f_{i}(t)f_{j}(t + \tau_{1})f_{k}(t + \tau_{2})] = R^{(3)}_{ijk}(\tau_{1}, \tau_{2}) \ for \ i, j, k=1, 2, \dots, m
\end{equation}
and can be represented through the following tensorial form as a multidimensional array
\begin{equation}
\begin{tikzpicture}[every node/.style={anchor=north east,fill=white,minimum width=1.4cm,minimum height=7mm}]
\matrix (mA) [draw,matrix of math nodes]
{
R^{(3)}_{11m}(\tau_{1}, \tau_{2}) & R^{(3)}_{12m}(\tau_{1}, \tau_{2}) & \dots  & R^{(3)}_{1mm}(\tau_{1}, \tau_{2}) \\
R^{(3)}_{21m}(\tau_{1}, \tau_{2}) & R^{(3)}_{22m}(\tau_{1}, \tau_{2}) & \dots  & R^{(3)}_{1mm}(\tau_{1}, \tau_{2}) \\
\vdots                & \vdots                & \ddots & \vdots                \\
R^{(3)}_{m1m}(\tau_{1}, \tau_{2}) & R^{(3)}_{m2m}(\tau_{1}, \tau_{2}) & \dots  & R^{(3)}_{mmm}(\tau_{1}, \tau_{2}) \\
};

\matrix (mB) [draw,matrix of math nodes] at ($(mA.south west)+(4.0,1.33)$)
{
R^{(3)}_{112}(\tau_{1}, \tau_{2}) & R^{(3)}_{122}(\tau_{1}, \tau_{2}) & \dots  & R^{(3)}_{1m2}(\tau_{1}, \tau_{2}) \\
R^{(3)}_{212}(\tau_{1}, \tau_{2}) & R^{(3)}_{222}(\tau_{1}, \tau_{2}) & \dots  & R^{(3)}_{1m2}(\tau_{1}, \tau_{2}) \\
\vdots                & \vdots                & \ddots & \vdots                \\
R^{(3)}_{m12}(\tau_{1}, \tau_{2}) & R^{(3)}_{m22}(\tau_{1}, \tau_{2}) & \dots  & R^{(3)}_{mm2}(\tau_{1}, \tau_{2}) \\
};

\matrix (mC) [draw,matrix of math nodes] at ($(mB.south west)+(7.5,2.5)$)
{
R^{(3)}_{111}(\tau_{1}, \tau_{2}) & R^{(3)}_{121}(\tau_{1}, \tau_{2}) & \dots  & R^{(3)}_{1m1}(\tau_{1}, \tau_{2}) \\
R^{(3)}_{211}(\tau_{1}, \tau_{2}) & R^{(3)}_{221}(\tau_{1}, \tau_{2}) & \dots  & R^{(3)}_{1m1}(\tau_{1}, \tau_{2}) \\
\vdots                & \vdots                & \ddots & \vdots                \\
R^{(3)}_{m11}(\tau_{1}, \tau_{2}) & R^{(3)}_{m21}(\tau_{1}, \tau_{2}) & \dots  & R^{(3)}_{mm1}(\tau_{1}, \tau_{2}) \\
};
\draw[dashed](mA.north east)--(mC.north east);
\draw[dashed](mA.north west)--(mC.north west);
\draw[dashed](mA.south east)--(mC.south east);
\end{tikzpicture}
\end{equation}

For real-valued, third-order stationary stochastic vector processes, the following second-order symmetry relationships hold,
\begin{equation}
\begin{aligned}
	& R_{ij}(\tau) = R_{ij}(-\tau), \ i, j = 1, 2, \dots, m \ \\
	& R_{ij}(\tau) = R_{ji}(\tau), \ i, j = 1, 2, \dots, m \ \\
\end{aligned}
\end{equation}
and the following third-order symmetry conditions hold,
\begin{equation}
\begin{aligned}
	& R^{(3)}_{ijk}(\tau_{1}, \tau_{2}) = R^{(3)}_{ijk}(\tau_{2}, \tau_{1}), \ i, j, k = 1, 2, \dots, m \\
	& R^{(3)}_{ijk}(\tau_{1}, \tau_{2}) = R^{(3)}_{ijk}(-\tau_{1}, -\tau_{2}), \ i, j, k = 1, 2, \dots, m \\
	& R^{(3)}_{ijk}(\tau_{1}, \tau_{2}) = R^{(3)}_{ijk}(-\tau_{1}, -\tau_{1} -\tau_{2}), \ i, j, k = 1, 2, \dots, m \\
	& R^{(3)}_{ijk}(\tau_{1}, \tau_{2}) = R^{(3)}_{ikj}(\tau_{1}, \tau_{2}) = R^{(3)}_{jik}(\tau_{1}, \tau_{2}) = R^{(3)}_{jki}(\tau_{1}, \tau_{2}) \\
	& = R^{(3)}_{kij}(\tau_{1}, \tau_{2}) = R^{(3)}_{kji}(\tau_{1}, \tau_{2}), \ i, j, k = 1, 2, \dots, m \\
\end{aligned}
\end{equation}

\subsection{Spectral Properties of Stochastic Vector Processes}

The ensemble properties of the stochastic vector properties relate to their spectral properties through the following Wiener-Khintchine transformations
\begin{equation}
    S_{jk}(\omega) = \frac{1}{2\pi}\int_{-\infty}^{\infty}R_{jk}(\tau)e^{-\iota \omega \tau}d\tau \ \text{for} \ j,k=1, 2, \dots, m
\end{equation}
\begin{equation}
    R_{ij}(\tau) = \int_{-\infty}^{\infty}S_{ij}(\omega)e^{\iota \omega \tau}d\tau \ \text{for} \ i,j=1, 2, \dots, m
\end{equation}
\begin{equation}
\begin{aligned}
    B_{ijk}(\omega_{1}, \omega_{2}) & = \frac{1}{(2\pi)^{2}}\int_{-\infty}^{\infty}R^{(3)}_{ijk}(\tau_{1}, \tau_{2})e^{-\iota (\omega_{1} \tau_{1} + \omega_{2}\tau_{2})}d\tau_{1}d\tau_{2}\
     \text{for} \ i,j,k=1, 2, \dots, m
\end{aligned}
\end{equation}
\begin{equation}
\begin{aligned}
    R^{(3)}_{ijk}(\tau_{1}, \tau_{2}) & = \int_{-\infty}^{\infty}B_{ijk}(\omega_{1}, \omega_{2})e^{\iota (\omega_{1} \tau_{1} + \omega_{2}\tau_{2})}d\omega_{1}d\omega_{2}\
    \text{for} \ i,j,k=1, 2, \dots, m
\end{aligned}
\end{equation}
where $B_{ijk}(\omega_{1}, \omega_{2})$ is the cross-bispectral density.

The second-order cross spectral density can be represented by following cross-spectral density matrix having terms given above:
\begin{equation}
S(\omega) = 
\begin{bmatrix}
    S_{11}(\omega) & S_{12}(\omega) & S_{13}(\omega) & \dots  & S_{1m}(\omega) \\
    S_{21}(\omega) & S_{22}(\omega) & S_{23}(\omega) & \dots  & S_{2m}(\omega) \\
    S_{31}(\omega) & S_{32}(\omega) & S_{33}(\omega) & \dots  & S_{3m}(\omega) \\
    \vdots        & \vdots        & \vdots        & \ddots & \vdots        \\
    S_{m1}(\omega) & S_{m2}(\omega) & S_{m3}(\omega) & \dots  & S_{mm}(\omega) \\
\end{bmatrix}
\end{equation}
Similarly, the third-order cross-bispectral density is represented in tensorial form by the following multidimensional array
\begin{equation}
\begin{tikzpicture}[every node/.style={anchor=north east,fill=white,minimum width=1.4cm,minimum height=7mm}]
\matrix (mA) [draw,matrix of math nodes]
{
B_{11m}(\omega_{1}, \omega_{2}) & B_{12m}(\omega_{1}, \omega_{2}) & \dots  & B_{1mm}(\omega_{1}, \omega_{2}) \\
B_{21m}(\omega_{1}, \omega_{2}) & B_{22m}(\omega_{1}, \omega_{2}) & \dots  & B_{1mm}(\omega_{1}, \omega_{2}) \\
\vdots                & \vdots                & \ddots & \vdots                \\
B_{m1m}(\omega_{1}, \omega_{2}) & B_{m2m}(\omega_{1}, \omega_{2}) & \dots  & B_{mmm}(\omega_{1}, \omega_{2}) \\
};

\matrix (mB) [draw,matrix of math nodes] at ($(mA.south west)+(4.0,1.33)$)
{
B_{112}(\omega_{1}, \omega_{2}) & B_{122}(\omega_{1}, \omega_{2}) & \dots  & B_{1m2}(\omega_{1}, \omega_{2}) \\
B_{212}(\omega_{1}, \omega_{2}) & B_{222}(\omega_{1}, \omega_{2}) & \dots  & B_{1m2}(\omega_{1}, \omega_{2}) \\
\vdots                & \vdots                & \ddots & \vdots                \\
B_{m12}(\omega_{1}, \omega_{2}) & B_{m22}(\omega_{1}, \omega_{2}) & \dots  & B_{mm2}(\omega_{1}, \omega_{2}) \\
};

\matrix (mC) [draw,matrix of math nodes] at ($(mB.south west)+(7.5,2.5)$)
{
B_{111}(\omega_{1}, \omega_{2}) & B_{121}(\omega_{1}, \omega_{2}) & \dots  & B_{1m1}(\omega_{1}, \omega_{2}) \\
B_{211}(\omega_{1}, \omega_{2}) & B_{221}(\omega_{1}, \omega_{2}) & \dots  & B_{1m1}(\omega_{1}, \omega_{2}) \\
\vdots                & \vdots                & \ddots & \vdots                \\
B_{m11}(\omega_{1}, \omega_{2}) & B_{m21}(\omega_{1}, \omega_{2}) & \dots  & B_{mm1}(\omega_{1}, \omega_{2}) \\
};
\draw[dashed](mA.north east)--(mC.north east);
\draw[dashed](mA.north west)--(mC.north west);
\draw[dashed](mA.south east)--(mC.south east);
\end{tikzpicture}
\end{equation}

The second and third order cross spectral density functions are complex valued in general and the following symmetry conditions hold

\begin{equation}
\begin{aligned}
    & S_{jj}(\omega) = S_{jj}(-\omega), \ j = 1, 2, \dots m \\
    & S_{ij}(\omega) = S_{ij}^{*}(-\omega), \ i, j = 1, 2, \dots m; \ i \neq j\\
    & S_{ij}(\omega) = S_{ji}^{*}(\omega), \ i, j = 1, 2, \dots m; \ i \neq j
\end{aligned}
\end{equation}
    
\begin{equation}
\begin{aligned}
    & B_{jjj}(\omega_{1}, \omega_{2}) = B_{jjj}(\omega_{2}, \omega_{1})\\
    & B_{jkl}(\omega_{1}, \omega_{2}) = B_{jkl}^{*}(\omega_{2}, \omega_{1})
\end{aligned}
\end{equation}

\section{Simulation of stochastic vector processes by 2nd-order Spectral Representation method}

The formula for the simulation of ergodic stochastic vector processes by $2^{nd}$-order Spectral Representation method is proposed in \cite{Deodatis1996}. In this method, only $2^{nd}$-order spectral information in the form of cross-spectral density is used for simulating sample functions. First, the cross spectral density $S(\omega)$ is decomposed as follows
\begin{equation}
    S(\omega) = H(\omega)H^{T*}(\omega)
    \label{eqn:cholesky}
\end{equation}
where the superscript ${}^{T*}$ denotes the matrix conjugate transpose. This decomposition can be carried out in many ways including Cholesky decomposition or eigenvalue decomposition. Throughout this study, whenever such a decomposition of any type of cross spectral density is presented it is assumed to be done by eigenvalue decomposition. This preference is due to the numerical stability of the eigenvalue decomposition in different computational implementations. In this case, the eigenvalue decomposition is given explicitly by
\begin{equation}
    S = \Phi \Sigma \Phi^T
\end{equation}
and the matrix $H$ in Eq.\ \eqref{eqn:cholesky} is given by $H=\Phi \sqrt{\Sigma}$.

The decomposed cross-spectral density matrix takes the following form
\begin{equation}
H(\omega) = 
\begin{bmatrix}
    H_{11}(\omega) & H_{12}(\omega) & H_{13}(\omega) & \dots  & H_{1m}(\omega) \\
    H_{21}(\omega) & H_{22}(\omega) & H_{23}(\omega) & \dots  & H_{2m}(\omega) \\
    H_{31}(\omega) & H_{32}(\omega) & H_{33}(\omega) & \dots  & H_{3m}(\omega) \\
    \vdots        & \vdots        & \vdots        & \ddots & \vdots        \\
    H_{m1}(\omega) & H_{m2}(\omega) & H_{m3}(\omega) & \dots  & H_{mm}(\omega) 
\end{bmatrix}
\end{equation}
where the diagonal elements are real functions of $\omega$ and the off-diagonal elements could potentially be complex functions of $\omega$ and possesses the following symmetry conditions,
\begin{equation}
    \begin{aligned}
    & H_{jj}(\omega) = H_{jj}(-\omega); \ j=1,2, \dots, m\\
    & H_{jk}(\omega) = H_{jk}^{*}(-\omega); \ j,k=1,2, \dots, m; \ j\neq k
    \end{aligned}
\end{equation}
The off-diagonal elements can be written in the polar form as
\begin{equation}
    \begin{aligned}
    & H_{jk}(\omega) = |H_{jk}(\omega)|e^{\iota\theta_{jk}(\omega)}; \ j,k=1,2, \dots, n; j\neq k
    \end{aligned}
\end{equation}
where
\begin{equation}
    \begin{aligned}
    & \theta_{jk}(\omega) = \arctan \left\{ \frac{\Im[H_{jk}(\omega)]}{\Re[H_{jk}(\omega)]} \right\} 
    \end{aligned}
\end{equation}

The stochastic process $f_{j}(t); j= 1,2, \dots n$ can be simulated by the following summation
\begin{equation}
    \begin{aligned}
    f_{j}(t) = \sum_{l=1}^{m}\sum_{k=0}^{N-1} |H_{jl}(\omega_{lk})|\sqrt{\Delta\omega}cos(\omega_{lk}t - \theta_{jl}(\omega_{lk}) + \phi_{lk}); \ j=1,2 \dots, m
    \end{aligned}
\end{equation}
where
\begin{equation}
    \begin{aligned}
    & \omega_{lk} = k\Delta\omega + \frac{l}{m}\Delta\omega\\
    & \theta_{jl}(\omega_{lk}) = \arctan \left\{ \frac{\Im[H_{jl}(\omega_{lk})]}{\Re[H_{jl}(\omega_{lk})]} \right\}
    \end{aligned}
\end{equation}
and $\phi_{lk}$ are indepedent random phase angles uniformly distribution on the interval $[0,2\pi]$.


The sample functions simulated by the above simulation formula are ergodic in first and second-order correlation properties, details of which can be found in \cite{Deodatis1996}.

\section{Simulation of ergodic $3^{rd}$-order random processes by Spectral Representation method}
\label{S:ergodicity}

Let us briefly return to the case of univariate stochastic processes to introduce a slight variation to the $3^{rd}$-order SRM presented in \cite{Shields2017} to generate ergodic sample functions. We then provide proofs of ergodicity and return to the multi-variate case in Section \ref{sec:sim_mv_SRm}.

The formula presented in \cite{Shields2017} for the simulation of $3^{rd}$-order stochastic processes with asymmetric non-linear wave interactions defined through a specified power spectrum $S(\omega)$ and bispectrum $B(\omega_i, \omega_j)$, is given by
\begin{equation}
\begin{aligned}
    & f(t) = 2\sum_{k=0}^{N-1} \sqrt{S_{p}(\omega_{k})\Delta\omega}\cos(\omega_{k}t + \phi_{k})\\
    & + 2\sum_{k=0}^{N-1}\sum_{i + j = k}^{i \geq j > 0}\frac{|B(\omega_{i}, \omega_{j})|}{\sqrt{S_{p}(\omega_{i})S_{p}(\omega_{j})}} \Delta \omega \cos(\omega_{i}t + \omega_{j}t + \phi_{i} + \phi_{j})
\end{aligned}
\label{eqn:3rd-order_SRm}
\end{equation}
where
\begin{equation}
\begin{aligned}
    S_{p}(\omega_{k}) = S(\omega_{k})\Big(1 - \sum_{i + j = k}^{i \geq j > 0}\frac{|B(\omega_{i}, \omega_{j})|^{2} \Delta \omega^{2}}{S_{p}(\omega_{i})S_{p}(\omega_{j})S(\omega_{i} + \omega_{j})} \Big),
\label{eqn:univaraite_pure_power_spectrum}
\end{aligned}
\end{equation}
\begin{equation}
    \omega_{k} = k \Delta \omega,
\end{equation}
and $\phi_i$ are independent random phase angles uniformly distributed on the interval $[0,2\pi]$.

Although the samples functions simulated by the above formula are proven to satisfy both the prescribed power spectrum and bispectrum in ensemble, they are not ergodic. But simply modifying the frequency indexing increment in the following way makes the sample functions ergodic up to the third order. 
\begin{equation}
\begin{aligned}
	&\omega_{k} = \left(k + \frac{1}{N}\right)\Delta\omega
\end{aligned}
\end{equation}
With the new increments, the sample functions are periodic, with period
\begin{equation}
\begin{aligned}
	T_{0} = &\frac{2N\pi}{\Delta\omega}
\label{eqn:uni_variate_ergodicity_time}
\end{aligned}
\end{equation}
and the sample functions are ergodic up to third order when simulated up to time $T = T_{0}$. The ergodicity proofs follow.




\subsection{Ergodicity in mean-value}

Ergodicity in the mean-value requires that:
\begin{equation}
    \langle f(t) \rangle_{T} = \mathbb{E}[f(t)]
\end{equation}
where $\langle \cdot \rangle_{T}$ denotes a time averaging over a time interval $T$. Ergodicity in the mean value is proven by substituting the simulation formula in Eq.\ \eqref{eqn:3rd-order_SRm} into the temporal averaging as
\begin{equation}
\begin{aligned}
	& \langle f(t) \rangle_{T} = \frac{1}{T}\int_{0}^{T}f(t)dt\\
	& =\frac{2}{T}\int_{0}^{T}\sum_{k=0}^{N-1} \sqrt{S_{p}(\omega_{k})\Delta \omega }\cos(\omega_{k}t + \phi_{k})\\
	& + \sum_{k=0}^{N-1}\sum_{i + j = k}^{i \geq j > 0}\frac{|B(\omega_{i}, \omega_{j})|}{\sqrt{S_{p}(\omega_{i})S_{p}(\omega_{j})}}\Delta \omega \cos(\omega_{i}t + \omega_{j}t + \phi_{i} + \phi_{j})dt\\
\end{aligned}
\end{equation}

The integrand in the above equation is periodic in $t$ with period equal to $T_{0}$ as given in \eqref{eqn:uni_variate_ergodicity_time}. It is therefore obvious that
\begin{equation}
\begin{aligned}
	&\langle f(t) \rangle_{T} = \mathbb{E}[f(t)] = 0 \ \text{when} \ T = T_{0}\\
\end{aligned}
\end{equation}

The above equation shows that the temporal average is equal to the ensemble value when the length of the sample function is equal to $T_{0}$. Hence, the sample functions are ergodic in the mean.

\subsection{Ergodicity in correlation}
Ergodicity in correlation requires that:
\begin{equation}
    \langle f(t) f(t + \tau) \rangle_{T} = E[f(t)f(t+\tau)] = R(\tau)
\end{equation}

The temporal average of the product of any two time instances of a sample function can be computed by, again substituting the simulation formal for $f(t)$ as follows
\begin{equation}
\begin{aligned}
	& \langle f(t) f(t + \tau) \rangle_{T} = \frac{1}{T}\int_{0}^{T}f(t)f(t + \tau)dt\\
	& =\frac{4}{T}\int_{0}^{T} \Big( \sum_{k_{1}=0}^{N-1} \sqrt{S_{p}(\omega_{k_{1}})\Delta \omega}\cos(\omega_{k_{1}}t + \phi_{k_{1}})\\
	& + \sum_{k_{1}=0}^{N-1}\sum_{i_{1} + j_{1} = k_{1}}^{i_{1} \geq j_{1} > 0}\frac{|B(\omega_{i_{1}}, \omega_{j_{1}})|}{\sqrt{S_{p}(\omega_{i_{1}})S_{p}(\omega_{j_{1}})}}\Delta \omega \cos((\omega_{i_{1}} + \omega_{j_{1}})t + \phi_{i_{1}} + \phi_{j_{1}}) \Big)\\
	& \Big( \sum_{k_{2}=0}^{N-1} \sqrt{S_{p}(\omega_{k_{2}})\Delta \omega}\cos(\omega_{k_{2}}(t + \tau) + \phi_{k_{2}}) \\
	& + \sum_{k_{2}=0}^{N-1}\sum_{i_{2} + j_{2} = k_{2}}^{i_{2} \geq j_{2} > 0}\frac{|B(\omega_{i_{2}}, \omega_{j_{2}})|}{\sqrt{S_{p}(\omega_{i_{2}})S_{p}(\omega_{j_{2}})}}\Delta \omega \cos((\omega_{i_{2}} + \omega_{j_{2}})(t + \tau) + \phi_{i_{2}} + \phi_{j_{2}}) \Big)dt
\end{aligned}
\end{equation}
Expanding the integrand yields:
\begin{equation}
\begin{aligned}
	& \langle f(t) f(t + \tau) \rangle_{T} = \frac{1}{T}\int_{0}^{T}f(t)f(t + \tau)dt\\
	& = \frac{4}{T}\int_{0}^{T} \Big( \sum_{k_{1}=0}^{N-1}\sum_{k_{2}=0}^{N-1} \sqrt{S_{p}(\omega_{k_{1}})\Delta \omega}\cos(\omega_{k_{1}}t + \phi_{k_{1}}) \sqrt{S_{p}(\omega_{k_{2}})\Delta \omega}\cos(\omega_{k_{2}}(t + \tau) + \phi_{k_{2}}) \\
	& + \sum_{k_{2}=0}^{N-1}\sum_{k_{1}=0}^{N-1}\sum_{i_{1} + j_{1} = k_{1}}^{i_{1} \geq j_{1} > 0}\sqrt{S_{p}(\omega_{k_{2}})\Delta \omega}\cos(\omega_{k_{2}}(t + \tau) + \phi_{k_{2}}) \frac{|B(\omega_{i_{1}}, \omega_{j_{1}})|}{\sqrt{S_{p}(\omega_{i_{1}})S_{p}(\omega_{j_{1}})}}\Delta \omega \cos((\omega_{i_{1}} + \omega_{j_{1}})t + \phi_{i_{1}} + \phi_{j_{1}})\\
	& + \sum_{k_{1}=0}^{N-1}\sum_{k_{2}=0}^{N-1}\sum_{i_{2} + j_{2} = k_{2}}^{i_{2} \geq j_{2} > 0}\sqrt{S_{p}(\omega_{k_{1}})\Delta \omega}\cos(\omega_{k_{1}}t + \phi_{k_{1}}) \frac{|B(\omega_{i_{2}}, \omega_{j_{2}})|}{\sqrt{S_{p}(\omega_{i_{2}})S_{p}(\omega_{j_{2}})}}\Delta \omega \cos((\omega_{i_{2}} + \omega_{j_{2}})(t + \tau) + \phi_{i_{2}} + \phi_{j_{2}})\\
	+ & \sum_{k_{1}=0}^{N-1}\sum_{k_{2}=0}^{N-1}\sum_{i_{1} + j_{1} = k_{1}}^{i_{1} \geq j_{1} > 0}\sum_{i_{2} + j_{2} = k_{2}}^{i_{2} \geq j_{2} > 0}\frac{|B(\omega_{i_{1}}, \omega_{j_{1}})|}{\sqrt{S_{p}(\omega_{i_{1}})S_{p}(\omega_{j_{1}})}}\Delta \omega \cos((\omega_{i_{1}} + \omega_{j_{1}})t + \phi_{i_{1}} + \phi_{j_{1}})\\
	& \frac{|B(\omega_{i_{2}}, \omega_{j_{2}})|}{\sqrt{S_{p}(\omega_{i_{2}})S_{p}(\omega_{j_{2}})}}\Delta \omega \cos((\omega_{i_{2}} + \omega_{j_{2}})(t + \tau) + \phi_{i_{2}} + \phi_{j_{2}}) \Big)dt
\end{aligned}
\end{equation}

The integrand of this equation has four terms. The first term can be simplified as follows:
\begin{equation}
\begin{aligned}
	&\frac{4}{T}\int_{0}^{T} \sum_{k_{1}=0}^{N-1}\sum_{k_{2}=0}^{N-1} \sqrt{S_{p}(\omega_{k_{1}})\Delta \omega}\cos(\omega_{k_{1}}t + \phi_{k_{1}}) \sqrt{S_{p}(\omega_{k_{2}})\Delta \omega}\cos(\omega_{k_{2}}(t + \tau) + \phi_{k_{2}}) \\
	& = \frac{4}{T}\int_{0}^{T} \sum_{k_{1}=0}^{N-1}\sum_{k_{2}=0}^{N-1} \sqrt{S_{p}(\omega_{k_{1}})S_{p}(\omega_{k_{2}})}\Delta \omega \cos(\omega_{k_{1}}t + \phi_{k_{1}})\cos(\omega_{k_{2}}(t + \tau) + \phi_{k_{2}}) \\
	& = \frac{2}{T}\int_{0}^{T} \sum_{k_{1}=0}^{N-1}\sum_{k_{2}=0}^{N-1} \sqrt{S_{p}(\omega_{k_{1}})S_{p}(\omega_{k_{2}})}\Delta \omega \Big( \cos(\omega_{k_{1}}t + \omega_{k_{2}}(t + \tau) +  \phi_{k_{1}} + \phi_{k_{2}}) \\
	& + \cos(\omega_{k_{1}}t - \omega_{k_{2}}(t + \tau) + \phi_{k_{1}} - \phi_{k_{2}})\Big) \\
	& = \frac{2}{T}\int_{0}^{T} \sum_{k_{1}=0}^{N-1}\sum_{k_{2}=0}^{N-1} \sqrt{S_{p}(\omega_{k_{1}})S_{p}(\omega_{k_{2}})} \Delta \omega \Big( \cos((\omega_{k_{1}} + \omega_{k_{2}})t + \omega_{k_{2}}\tau +  \phi_{k_{1}} + \phi_{k_{2}}) \\
	& + \cos((\omega_{k_{1}} - \omega_{k_{2}})t - \omega_{k_{2}}\tau + \phi_{k_{1}} - \phi_{k_{2}})\Big) 
\end{aligned}
\end{equation}
It is straightforward to see that
\begin{equation}
\begin{aligned}
	&\int_{0}^{T} \cos((\omega_{k_{1}} + \omega_{k_{2}})t)dt = 0; \ \forall  k_{1}, k_{2} \\
	&\int_{0}^{T} \cos((\omega_{k_{1}} - \omega_{k_{2}})t)dt = 0; \ \forall  k_{1}, k_{2}, \ k_{1} \neq \ k_{2} \\
\end{aligned}
\end{equation}
so the expression simplifies to
\begin{equation}
\begin{aligned}
	&\frac{2}{T}\int_{0}^{T} \sum_{k=0}^{N-1} S_{p}(\omega_{k}) \Delta \omega \cos(- \omega_{k}\tau) dt \\
	& = 2\sum_{k=0}^{N-1} S_{p}(\omega_{k}) \Delta \omega \cos(- \omega_{k}\tau) \\
	& = 2\sum_{k=0}^{N-1} S_{p}(\omega_{k}) \Delta \omega \cos( \omega_{k}\tau) \\
\label{eqn:second_order_p}
\end{aligned}
\end{equation}

The second term reduces as follows:
\begin{equation}
\begin{aligned}
	& \frac{4}{T}\int_{0}^{T} \sum_{k_{2}=0}^{N-1}\sum_{k_{1}=0}^{N-1}\sum_{i_{1} + j_{1} = k_{1}}^{i_{1} \geq j_{1} > 0}\sqrt{S_{p}(\omega_{k_{2}})\Delta \omega}\cos(\omega_{k_{2}}(t + \tau) + \phi_{k_{2}})\\
	& \frac{|B(\omega_{i_{1}}, \omega_{j_{1}})|}{\sqrt{S_{p}(\omega_{i_{1}})S_{p}(\omega_{j_{1}})}}\Delta \omega \cos((\omega_{i_{1}} + \omega_{j_{1}})t + \phi_{i_{1}} + \phi_{j_{1}})dt\\
	& = \frac{4}{T}\int_{0}^{T} \sum_{k_{2}=0}^{N-1}\sum_{k_{1}=0}^{N-1}\sum_{i_{1} + j_{1} = k_{1}}^{i_{1} \geq j_{1} > 0}\sqrt{S_{p}(\omega_{k_{2}})\Delta \omega}\frac{|B(\omega_{i_{1}}, \omega_{j_{1}})|}{\sqrt{S_{p}(\omega_{i_{1}})S_{p}(\omega_{j_{1}})}}\Delta \omega\\
	& \cos(\omega_{k_{2}}(t + \tau) + \phi_{k_{2}})\cos((\omega_{i_{1}} + \omega_{j_{1}})t + \phi_{i_{1}} + \phi_{j_{1}})dt\\
	& = \frac{2}{T}\int_{0}^{T} \sum_{k_{2}=0}^{N-1}\sum_{k_{1}=0}^{N-1}\sum_{i_{1} + j_{1} = k_{1}}^{i_{1} \geq j_{1} > 0}\sqrt{S_{p}(\omega_{k_{2}})\Delta \omega}\frac{|B(\omega_{i_{1}}, \omega_{j_{1}})|}{\sqrt{S_{p}(\omega_{i_{1}})S_{p}(\omega_{j_{1}})}}\Delta \omega\\
	& \Big( \cos(\omega_{k_{2}}(t + \tau) + (\omega_{i_{1}} + \omega_{j_{1}})t + \phi_{k_{2}} + \phi_{i_{1}} + \phi_{j_{1}}) \\
	& + \cos(\omega_{k_{2}}(t + \tau) + (\omega_{i_{1}} + \omega_{j_{1}})t + \phi_{k_{2}} + \phi_{i_{1}} + \phi_{j_{1}}) \Big)dt\\
	& = \frac{2}{T}\int_{0}^{T} \sum_{k_{2}=0}^{N-1}\sum_{k_{1}=0}^{N-1}\sum_{i_{1} + j_{1} = k_{1}}^{i_{1} \geq j_{1} > 0}\sqrt{S_{p}(\omega_{k_{2}})\Delta \omega}\frac{|B(\omega_{i_{1}}, \omega_{j_{1}})|}{\sqrt{S_{p}(\omega_{i_{1}})S_{p}(\omega_{j_{1}})}}\Delta \omega\\
	& \Big( \cos((\omega_{i_{1}} + \omega_{j_{1}} + \omega_{k_{2}})t  + \omega_{k_{2}}\tau + \phi_{k_{2}} + \phi_{i_{1}} + \phi_{j_{1}}) \\
	& + \cos((\omega_{k_{2}} - \omega_{i_{1}} - \omega_{j_{1}})t + \omega_{k_{2}}\tau + \phi_{k_{2}} - \phi_{i_{1}} - \phi_{j_{1}}) \Big)dt\\
	=0
\end{aligned}
\end{equation}
because
\begin{equation}
\begin{aligned}
	& \int_{0}^{T} \cos((\omega_{i_{1}} + \omega_{j_{1}} + \omega_{k_{2}})t) dt = 0; \ \forall i_{1}, j_{1} \ and \ k_{2}\\
	& \int_{0}^{T} \cos((\omega_{k_{2}} - \omega_{i_{1}} - \omega_{j_{1}})t) dt = 0; \ \forall i_{1}, j_{1} \ and \ k_{2}
\end{aligned}
\end{equation}
since
\begin{equation}
    \omega_{k_{2}} - \omega_{i_{1}} - \omega_{j_{1}} = \left(k_{2} - i_{1} - j_{1} + 1 + \frac{k_{2} - i_{1} - j_{1}}{N}\right).
\end{equation}
The third term of the integrand similarly reduces to zero.

The final term is reduces as follows:
\begin{equation}
\begin{aligned}
	& \frac{4}{T}\int_{0}^{T}\sum_{k_{1}=0}^{N-1}\sum_{k_{2}=0}^{N-1}\sum_{i_{1} + j_{1} = k_{1}}^{i_{1} \geq j_{1} > 0}\sum_{i_{2} + j_{2} = k_{2}}^{i_{2} \geq j_{2} > 0}\frac{|B(\omega_{i_{1}}, \omega_{j_{1}})|}{\sqrt{S_{p}(\omega_{i_{1}})S_{p}(\omega_{j_{1}})}}\Delta \omega \cos((\omega_{i_{1}} + \omega_{j_{1}})t + \phi_{i_{1}} + \phi_{j_{1}})\\
	& \frac{|B(\omega_{i_{2}}, \omega_{j_{2}})|}{\sqrt{S_{p}(\omega_{i_{2}})S_{p}(\omega_{j_{2}})}}\Delta \omega \cos((\omega_{i_{2}} + \omega_{j_{2}})(t + \tau) + \phi_{i_{2}} + \phi_{j_{2}}) dt\\
	& = \frac{4}{T}\int_{0}^{T}\sum_{k_{1}=0}^{N-1}\sum_{k_{2}=0}^{N-1}\sum_{i_{1} + j_{1} = k_{1}}^{i_{1} \geq j_{1} > 0}\sum_{i_{2} + j_{2} = k_{2}}^{i_{2} \geq j_{2} > 0}\frac{|B(\omega_{i_{1}}, \omega_{j_{1}})|}{\sqrt{S_{p}(\omega_{i_{1}})S_{p}(\omega_{j_{1}})}} \frac{|B(\omega_{i_{2}}, \omega_{j_{2}})|}{\sqrt{S_{p}(\omega_{i_{2}})S_{p}(\omega_{j_{2}})}}{\Delta \omega}^{2}\\
	& \cos((\omega_{i_{1}} + \omega_{j_{1}})t + \phi_{i_{1}} + \phi_{j_{1}})\cos((\omega_{i_{2}} + \omega_{j_{2}})(t + \tau) + \phi_{i_{2}} + \phi_{j_{2}}) dt\\
	& = \frac{2}{T}\int_{0}^{T}\sum_{k_{1}=0}^{N-1}\sum_{k_{2}=0}^{N-1}\sum_{i_{1} + j_{1} = k_{1}}^{i_{1} \geq j_{1} > 0}\sum_{i_{2} + j_{2} = k_{2}}^{i_{2} \geq j_{2} > 0}\frac{|B(\omega_{i_{1}}, \omega_{j_{1}})|}{\sqrt{S_{p}(\omega_{i_{1}})S_{p}(\omega_{j_{1}})}} \frac{|B(\omega_{i_{2}}, \omega_{j_{2}})|}{\sqrt{S_{p}(\omega_{i_{2}})S_{p}(\omega_{j_{2}})}}{\Delta \omega}^{2}\\
	& \Big( \cos((\omega_{i_{1}} + \omega_{j_{1}} + \omega_{i_{2}} + \omega_{j_{2}})t + (\omega_{i_{2}} + \omega_{j_{2}})\tau + \phi_{i_{1}} + \phi_{j_{1}} + \phi_{i_{2}} + \phi_{j_{2}})\\
	& + \cos((\omega_{i_{2}} + \omega_{j_{2}} - \omega_{i_{1}} - \omega_{j_{1}})t + (\omega_{i_{2}} + \omega_{j_{2}})\tau - \phi_{i_{1}} - \phi_{j_{1}} + \phi_{i_{2}} + \phi_{j_{2}}) \Big) dt\\
\end{aligned}
\end{equation}
We observe that
\begin{equation}
\begin{aligned}
	& \int_{0}^{T} \cos((\omega_{i_{1}} + \omega_{j_{1}} + \omega_{i_{2}} + \omega_{j_{2}})t) dt = 0; \ \forall i_{1}, j_{1}, i_{2}, j_{2}\\
	& \int_{0}^{T} \cos((\omega_{i_{2}} + \omega_{j_{2}} - \omega_{i_{1}} - \omega_{j_{1}})t) dt = 0; \ \forall i_{1}, j_{1}, i_{2}, j_{2}, \ i_{1} \neq i_{2}, j_{1} \neq j_{2}\\
\end{aligned}
\end{equation}
and therefore this simplifies to 
\begin{equation}
\begin{aligned}
	& \frac{2}{T}\int_{0}^{T}\sum_{k=0}^{N-1}\sum_{i + j = k}^{i \geq j > 0}\frac{|B(\omega_{i}, \omega_{j})|^{2}}{S_{p}(\omega_{i})S_{p}(\omega_{j})}{\Delta \omega}^{2} \cos((\omega_{i} + \omega_{j})\tau) dt\\
\end{aligned}
\end{equation}
Recognizing that $\omega_{k} = \omega_{i} + \omega_{j}$ when $k = i + j$, leads to
\begin{equation}
\begin{aligned}
	& \frac{2}{T}\int_{0}^{T}\sum_{k=0}^{N-1}\sum_{i + j = k}^{i \geq j > 0}\frac{|B(\omega_{i}, \omega_{j})|^{2}}{S_{p}(\omega_{i})S_{p}(\omega_{j})}{\Delta \omega}^{2}\cos(\omega_{k}\tau) dt \\
	& = \frac{2}{T}\int_{0}^{T}\sum_{k=0}^{N-1}S_{I}(\omega_{k}) \Delta \omega \cos(\omega_{k}\tau) dt \\
	& = 2 \sum_{k=0}^{N-1}S_{I}(\omega_{k}) \Delta \omega \cos(\omega_{k}\tau) \\
\label{eqn:second_order_i}
\end{aligned}
\end{equation}
where $S_I(\omega)$ denotes the portion of the power spectrum arising from wave interactions with $S(\omega)=S_p(\omega)+S_I(\omega)$.

Combining Eqs.\ \eqref{eqn:second_order_p} and \eqref{eqn:second_order_i}, we see that
\begin{equation}
\begin{aligned}
	& \langle f(t) f(t + \tau) \rangle_{T} = 2 \sum_{k=0}^{N-1}S(\omega_{k}) \Delta \omega \cos(\omega_{k}\tau) \\
	& \langle f(t) f(t + \tau) \rangle_{T} = R_{2}(\tau) \\
\end{aligned}
\end{equation}
and therefore the sample functions are ergodic in the second-order correlation when simulated with length $T_0$.

\subsection{Ergodicity in third-order moment}

Ergodicity in the third-order moment requires that:
\begin{equation}
    \langle f(t) f(t + \tau_1) f(t+\tau_2) \rangle_{T} = E[f(t)f(t+\tau_1)f(t+\tau_2)] = R^{(3)}(\tau_1,\tau_2)
\end{equation}
Again, we start by taking the temporal average of the product of any three time instances of a sample function as follows
\begin{equation}
\begin{aligned}
	& \langle f(t) f(t + \tau_{1}) f(t + \tau_{2}) \rangle_{T} = \frac{1}{T}\int_{0}^{T}f(t)f(t + \tau_{1})f(t + \tau_{2})dt\\
	& =\frac{8}{T}\int_{0}^{T} \Big( \sum_{k_{1}=0}^{N-1} \sqrt{S_{p}(\omega_{k_{1}})\Delta \omega}\cos(\omega_{k_{1}}t + \phi_{k_{1}})\\
	& + \sum_{k_{1}=0}^{N-1}\sum_{i_{1} + j_{1} = k_{1}}^{i_{1} \geq j_{1} > 0}\frac{|B(\omega_{i_{1}}, \omega_{j_{1}})|}{\sqrt{S_{p}(\omega_{i_{1}})S_{p}(\omega_{j_{1}})}}\Delta \omega \cos((\omega_{i_{1}} + \omega_{j_{1}})t + \phi_{i_{1}} + \phi_{j_{1}}) \Big)\\
	& \Big( \sum_{k_{2}=0}^{N} \sqrt{S_{p}(\omega_{k_{2}})\Delta \omega}\cos(\omega_{k_{2}}(t + \tau_{1}) + \phi_{k_{2}}) \\
	& + \sum_{k_{2}=0}^{N-1}\sum_{i_{2} + j_{2} = k_{2}}^{i_{2} \geq j_{2} > 0}\frac{|B(\omega_{i_{2}}, \omega_{j_{2}})|}{\sqrt{S_{p}(\omega_{i_{2}})S_{p}(\omega_{j_{2}})}}\Delta \omega \cos((\omega_{i_{2}} + \omega_{j_{2}})(t + \tau_{2}) + \phi_{i_{2}} + \phi_{j_{2}}) \Big)\\
	& \Big( \sum_{k_{3}=0}^{N-1} \sqrt{S_{p}(\omega_{k_{3}})\Delta \omega}\cos(\omega_{k_{3}}(t + \tau_{2}) + \phi_{k_{3}}) \\
	& + \sum_{k_{3}=0}^{N-1}\sum_{i_{3} + j_{3} = k_{3}}^{i_{3} \geq j_{3} > 0}\frac{|B(\omega_{i_{3}}, \omega_{j_{3}})|}{\sqrt{S_{p}(\omega_{i_{3}})S_{p}(\omega_{j_{3}})}}\Delta \omega \cos((\omega_{i_{3}} + \omega_{j_{3}})(t + \tau_{2}) + \phi_{i_{3}} + \phi_{j_{3}}) \Big) dt
\end{aligned}
\end{equation}
Only analyzing the terms that do not integrate to zero, we have
\begin{equation}
\begin{aligned}
	& \langle f(t) f(t + \tau_{1}) f(t + \tau_{2}) \rangle_{T}  =\frac{24}{T}\int_{0}^{T} \sum_{k_{3}=0}^{N-1} \sqrt{S_{p}(\omega_{k_{3}})\Delta \omega}\cos(\omega_{k_{3}}(t + \tau_{2}) + \phi_{k_{3}})\\
	& \sum_{k_{2}=0}^{N-1} \sqrt{S_{p}(\omega_{k_{2}})\Delta \omega}\cos(\omega_{k_{2}}(t + \tau_{1}) + \phi_{k_{2}})
	\sum_{k_{1}=0}^{N-1}\sum_{i_{1} + j_{1} = k_{1}}^{i_{1} \geq j_{1} > 0}\frac{|B(\omega_{i_{1}}, \omega_{j_{1}})|}{\sqrt{S_{p}(\omega_{i_{1}})S_{p}(\omega_{j_{1}})}}\Delta \omega \cos((\omega_{i_{1}} + \omega_{j_{1}})t + \phi_{i_{1}} + \phi_{j_{1}}) dt
\end{aligned}
\end{equation}
When $k_{2} = i_{1} = i$, $k_{3} = j_{1} = j$ and $k_{1} = k$, this reduces as
\begin{equation}
\begin{aligned}
	& \frac{48}{T}\int_{0}^{T} \sum_{k=0}^{N-1}\sum_{i + j = k}^{i \geq j > 0} \sqrt{S_{p}(\omega_{j})}\cos(\omega_{j}(t + \tau_{2}) + \phi_{j})\sqrt{S_{p}(\omega_{i})}\cos(\omega_{i}(t + \tau_{1}) + \phi_{i})\\
	&\frac{|B(\omega_{i}, \omega_{j})|}{\sqrt{S_{p}(\omega_{i})S_{p}(\omega_{j})}}{\Delta \omega}^{2}\cos((\omega_{i} + \omega_{j})(t) + \phi_{i} + \phi_{j}) dt\\
	& =\frac{48}{T}\int_{0}^{T} \sum_{k=0}^{N-1}\sum_{i + j = k}^{i \geq j > 0} |B(\omega_{i}, \omega_{j})|{\Delta \omega}^{2}\cos(\omega_{j}(t + \tau_{2}) + \phi_{j}) \cos(\omega_{i}(t + \tau_{1}) + \phi_{i})\\
	& \cos((\omega_{i} + \omega_{j})t + \phi_{i} + \phi_{j}) dt\\
	& =\frac{24}{T}\int_{0}^{T} \sum_{k=0}^{N-1}\sum_{i + j = k}^{i \geq j > 0} |B(\omega_{i}, \omega_{j})|{\Delta \omega}^{2}\Big( \cos((\omega_{i} + \omega_{j})t + \omega_{j}\tau_{2} + \omega_{i}\tau_{1} + \phi_{i} + \phi_{j}) \\
	& +  \cos((\omega_{i} - \omega_{j})t + \omega_{i}\tau_{1} - \omega_{j}\tau_{2} + \phi_{i} - \phi_{j}) \Big) \\
	& \cos((\omega_{i} + \omega_{j})t + \phi_{i} + \phi_{j}) dt\\
	& =\frac{12}{T}\int_{0}^{T} \sum_{k=0}^{N-1}\sum_{i + j = k}^{i \geq j > 0} |B(\omega_{i}, \omega_{j})|{\Delta \omega}^{2}\Big( \cos(2(\omega_{i} + \omega_{j})t + \omega_{j}\tau_{2} + \omega_{i}\tau_{1} + 2(\phi_{i} + \phi_{j})) \\
	& + \cos(\omega_{j}\tau_{2} + \omega_{i}\tau_{1}) + \cos(2\omega_{i}t+ \omega_{i}\tau_{1} - \omega_{j}\tau_{2} + 2\phi_{i}) + \cos(2\omega_{j}t+ \omega_{j}\tau_{2} - \omega_{i}\tau_{1} + 2\phi_{j}) \Big) dt
\end{aligned}
\end{equation}
We have that
\begin{equation}
\begin{aligned}
	&\int_{0}^{T}\cos(2(\omega_{i} + \omega_{j})t) = 0\\
	&\int_{0}^{T}\cos(2\omega_{i}t) = 0\\
	&\int_{0}^{T}\cos(2\omega_{j}t) = 0\\
\end{aligned}
\end{equation}
and therefore
\begin{equation}
\begin{aligned}
	\langle f(t) f(t + \tau_{1}) f(t + \tau_{2}) \rangle_{T} = 
	& = 12 \sum_{k=0}^{N-1}\sum_{i + j = k}^{i \geq j > 0} |B(\omega_{i}, \omega_{j})|{\Delta \omega}^{2}\cos(\omega_{j}\tau_{2} + \omega_{i}\tau_{1})\\
	& = 6 \sum_{k=0}^{N-1}\sum_{i + j = k}^{i, j > 0} |B(\omega_{i}, \omega_{j})|{\Delta \omega}^{2}\cos(\omega_{j}\tau_{2} + \omega_{i}\tau_{1})\\
	& = R^{3}(\tau_{1}, \tau_{2})
\end{aligned}
\end{equation}
As a result, the proposed simulation formula produces samples that are ergodic in their third-order moments.


\section{Simulation of ergodic stochastic vector processes by 3rd-order Spectral Representation method}
\label{sec:sim_mv_SRm}

Now that we have established the equation for simulation of univariate ergodic stochastic processes, let us return to the multi-variate case. To simulate the 1D-mV stationary stochastic vector process $[f_{1}(t), f_{2}(t), ....... f_{m}(t)]^{T}$, the pure component of the 2nd-order cross spectral density $S_{p}(\omega)$ must be computed first. In case of a simple 1D-1V stationary stochastic processes, the computation is straightforward, see Eq.\  \eqref{eqn:univaraite_pure_power_spectrum}. However, the pure cross-spectral density for a stochastic vector process is not trivial, and requires us to resort to Einstein (tensor) notation. This notation can be leveraged in various programming language like Python which make computations via the Einstein notation very appealing.
\begin{equation}
\begin{aligned}
    & S^{p}_{ab}(\omega_{k}) = S_{ab}(\omega_{k}) - \sum_{i + j = k}^{i \geq j > 0} B_{aef}(\omega_{i}, \omega_{j}) B^{*}_{bgh}(\omega_{i}, \omega_{j}) G_{pe}(\omega_{i}) G_{pg}(\omega_{i}) G_{qf}(\omega_{j}) G_{qh}(\omega_{j}) {\Delta\omega}^{2} \\
\end{aligned}
\end{equation}
where the term $G(\omega)$ is the inverse of the decomposed pure cross-spectral density derived as follows. Similar to the $2^{nd}$-order expansion, the pure cross-spectral density can be decomposed using the eigenvalue decomposition as
\begin{equation}
\begin{aligned}
    & S^{(p)}(\omega) = H(\omega)H^{T*}(\omega)\\
\end{aligned}
\end{equation}
and having the following properties
\begin{equation}
\begin{aligned}
    & H_{jj}(\omega) = H_{jj}(-\omega)\\
    & H_{jk}(\omega) = H_{jk}^{*}(-\omega)\\
    & H_{jk}(\omega) = |H_{jk}(\omega)|e^{\iota \theta_{jk}(\omega)}\\
    & \theta_{jk}(\omega) = \tan^{-1}\Big(\frac{\Im[H_{jk}(\omega)]}{\Re[H_{jk}(\omega)]}\Big)\\
\end{aligned}
\end{equation}
We then define $G(\omega) = (H(\omega))^{-1}$, which again can be expressed in polar coordinates as:
\begin{equation}
\begin{aligned}
    & G_{jk}(\omega) = |G_{jk}(\omega)|e^{\iota \theta^I_{jk}(\omega)}\\
    & \theta^{I}_{jk}(\omega) = \tan^{-1}\Big(\frac{\Im[G_{jk}(\omega)]}{\Re[G_{jk}(\omega)]}\Big)\\
\end{aligned}
\end{equation}




The stochastic process $f_{a}(t); a = 1, 2, \dots m$ in this case can be simulated as follows:
\begin{equation}
\begin{aligned}
    f_{a}(t) &= 2 \sum_{k=0}^{N-1} \big[\sum_{l=1}^{m}|H_{al}(\omega_{lk})|\sqrt{\Delta \omega}\cos(\omega_{lk}t - \theta_{al}(\omega_{lk}) + \phi_{lk}) \\
    & + 2 \sum_{l=1}^{m}\sum_{n=1}^{m}\sum_{p=1}^{m}\sum_{q=1}^{m}\sum_{i + j = k}^{i \geq j\geq 0} |B_{aln}(\omega_{pi}, \omega_{qj})||G_{lp}(\omega_{pi})||G_{nq}(\omega_{qj})| \Delta \omega\\
    & \cos((\omega_{pi} + \omega_{qj})t - \beta_{aln}(\omega_{pi}, \omega_{qj}) - \theta^{I}_{lp}(\omega_{pi}) - \theta^{I}_{nq}(\omega_{qj}) + \phi_{pi} + \phi_{qj})\big]
\end{aligned}
\label{eqn:mv_3rd_order_SRm}
\end{equation}
where
\begin{equation}
    \begin{aligned}
    & \beta_{aln}(\omega_{pi}, \omega_{qj}) = \tan^{-1}\Big(\frac{\Im[B_{aln}(\omega_{pi}, \omega_{qj})]}{\Re[B_{aln}(\omega_{pi}, \omega_{qj})]}\Big)
    \end{aligned}
\end{equation}
is the biphase, and 
\begin{equation}
    \begin{aligned}
    & \omega_{lk} = k\Delta\omega + \frac{l}{2m}\Delta\omega + \frac{1}{N}\Delta\omega
    \end{aligned}
    \label{eqn:multi-index}
\end{equation}

The stochastic vector processes simulated using Eq.\ \eqref{eqn:mv_3rd_order_SRm} satisfy both ensemble and ergodic properties of the vector process up to the third order. Proofs are presented in \ref{A:1} and \ref{A:2} respectively.

\subsection{Orthogonal increments involved in the simulation formula}

For completeness, the above simulation formula has been derived in the same manner as the univariate equation in \cite{Shields2017} and the multi-dimensional formula in Part I using a third-order orthogonal increment in the Cramer spectral representation. The orthogonal increment employed in the derivation of $3^{rd}$-order stochastic vector process is given by
\begin{equation}
\begin{aligned}
    & du_{a}(\omega_{k}) = 2\sum_{b=1}^{m} H_{ab}(\omega_{k}) cos(\phi_bk) +\\
    & 2\sum_{i+j=k}^{i \geq j \geq 0} \sum_{b=1}^{m} \sum_{c=1}^{m} \sum_{p=1}^{m} \sum_{q=1}^{m} B_{abc}(\omega_{i}, \omega_{j}) G_{bp}(\omega_{i})G_{cq}(\omega_{j})cos(\phi_{pi} + \phi_{qj}) \\
\end{aligned}
\end{equation}
It can be further verified that the orthogonal increment satisfies all orthogonality conditions up to third order presented in Part 1 \cite{VandanapuI}.

\section{Simulation by Fast Fourier Transform}
\label{S:FFT_implementation}

The simulation formula presented above for stochastic vector processes is computationally expensive, but an be accelerated with the Fast Fourier Transform (FFT). This section presents how the simulation formula can be expressed in a form suitable for FFT. The simulation formula is given in conventional form in Eq.\ \eqref{eqn:mv_3rd_order_SRm}. Applying Euler's formula, $e^{\iota\phi} = \cos(\phi) + \iota \sin(\phi)$, such that $\mathbb{R}[e^{\iota\phi}] = \cos(\phi)$, the simulation formula simplifies to 
\begin{equation}
\begin{aligned}
    f_{a}(t) &= \mathbb{R} \big[ 2 \sum_{k=0}^{N-1} [\sum_{l=1}^{m}|H_{al}(\omega_{lk})|\sqrt{\Delta \omega} e^{\iota(\omega_{lk}t - \theta_{al}(\omega_{lk}) + \phi_{lk})} \\
    & + \sum_{l=1}^{m}\sum_{n=1}^{m}\sum_{p=1}^{m}\sum_{q=1}^{m}\sum_{i + j = k}^{i \geq j\geq 0} |B_{aln}(\omega_{pi}, \omega_{qj})||G_{lp}(\omega_{pi})||G_{nq}(\omega_{qj})| \Delta \omega\\
    & e^{\iota (\omega_{pi} + \omega_{qj})t - \beta_{aln}(\omega_{pi}, \omega_{qj}) - \theta^{I}_{lp}(\omega_{pi}) - \theta^{I}_{nq}(\omega_{qj}) + \phi_{pi} + \phi_{qj}} ] \big]
\end{aligned}
\end{equation}
Discretizing the time domain using $t_r=r\Delta t$ and the frequency domain using the multi-indexed frequency in Eq.\ \eqref{eqn:multi-index} yields
\begin{equation}
\begin{aligned}
    f_{a}(r \Delta t) &= \mathbb{R} \big[ 2 \sum_{k=0}^{N-1} [ \sum_{l=1}^{m} |H_{al}(\omega_{lk})| \sqrt{\Delta \omega} e^{\iota (\frac{lr}{2m} + \frac{1}{N}) \Delta \omega \Delta t} e^{\iota(- \theta_{al}(\omega_{lk}) + \phi_{lk})}\\
    & + \sum_{l=1}^{m} \sum_{n=1}^{m} \sum_{p=1}^{m} \sum_{q=1}^{m} \sum_{i + j = k}^{i \geq j\geq 0} |B_{aln}(\omega_{i}, \omega_{j})||G_{lp}(\omega_{i})||G_{nq}(\omega_{j})| \Delta \omega\\
    & e^{\iota (\frac{pr}{2m} + \frac{qr}{2m} + \frac{2}{N}) \Delta \omega \Delta t} e^{\iota (- \beta_{aln}(\omega_{pi}, \omega_{qj}) - \theta^{I}_{lp}(\omega_{pi}) - \theta^{I}_{nq}(\omega_{qj}) + \phi_{pi} + \phi_{qj}}] e^{\iota kr \Delta \omega \Delta t} \big]
\end{aligned}
\end{equation}
Expressing this equation in terms of the standard FFT implementation, we have
\begin{equation}
\begin{aligned}
    f_{a}(r \Delta t) &= \mathbb{R} \big[ \sum_{k=0}^{N-1} C_{k} e^{\iota kr \Delta \omega \Delta t} \big]
\end{aligned}
\label{eqn:FFT1}
\end{equation}
where $C_{k}$ is given by
\begin{equation}
\begin{aligned}
    C_{k} &=  2 [ \sum_{l=1}^{m} |H_{al}(\omega_{lk})| \sqrt{\Delta \omega} e^{\iota(- \theta_{al}(\omega_{lk}) + \phi_{lk})} \\
    & + \sum_{i + j = k}^{i \geq j\geq 0} \sum_{l=1}^{m} \sum_{n=1}^{m} \sum_{p=1}^{m} \sum_{q=1}^{m} |B_{aln}(\omega_{pi}, \omega_{qj})||G_{lp}(\omega_{pi})||G_{nq}(\omega_{qj})| \Delta \omega\\
    & e^{\iota (- \beta_{aln}(\omega_{pi}, \omega_{qj}) - \theta^{I}_{lp}(\omega_{pi}) - \theta^{I}_{nq}(\omega_{qj}) + \phi_{pi} + \phi_{qj}} ]
\end{aligned}
\label{eqn:FFT2}
\end{equation}
where we have the following conditions to ensure the ergodicity and avoid aliasing
\begin{equation}
    mm\Delta t=T_0=m\frac{2\pi}{\Delta \omega}
\end{equation}
and 
\begin{equation}
    \Delta t = \frac{2\pi}{m\Delta \omega}.
\end{equation}
Simulation by FFT using Eq.\ \eqref{eqn:FFT1} and \eqref{eqn:FFT2} saves considerable computational expense while retaining the desired ensemble and ergoicity properties of the sample functions.

\section{Numerical Example}
\label{S:numerical_examples}

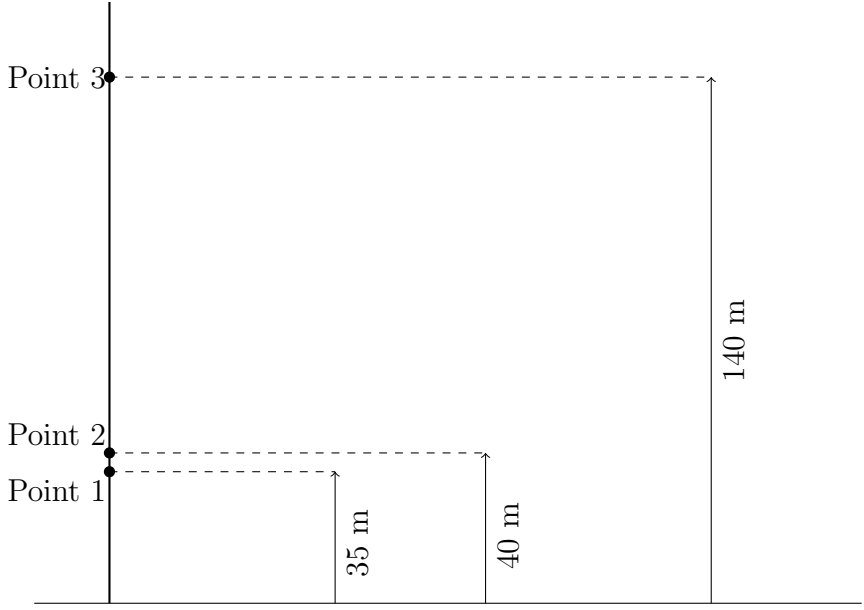
\begin{figure}[ht!]
\begin{tikzpicture}[scale = 1]
    \draw (-1, 0) -- (10, 0);
    \draw [thick] (0, 0) -- (0, 8);
    \draw [dashed] (0,1.75) -- (3,1.75);
    \draw [dashed] (0,2) -- (5,2);
    \draw [dashed] (0,7) -- (8,7);
    \draw [->] (3, 0) -- (3, 1.75);
    \draw [->] (5, 0) -- (5, 2);
    \draw [->] (8, 0) -- (8, 7);
    \tkzDefPoint (0,1.75){A}
    \tkzDefPoint (0,2){B}
    \tkzDefPoint (0,7){C}
    \tkzLabelPoint[left](A){}
    \tkzLabelPoint[left](B){}
    \tkzLabelPoint[left](C){}
    \node at (A) [circle,fill,inner sep=1.5pt]{};
    \node at (B) [circle,fill,inner sep=1.5pt]{};
    \node at (C) [circle,fill,inner sep=1.5pt]{};
    \node[rotate=90] at (3.3,0.8) {35 m};
    \node[rotate=90] at (5.3,0.9) {40 m};
    \node[rotate=90] at (8.3,3.5) {140 m};
    \node at (-0.7,1.5) {Point 1};
    \node at (-0.7,2.25) {Point 2};
    \node at (-0.7,7) {Point 3};
\end{tikzpicture}
\caption{Configuration of the wind velocity points along a vertical wind profile.}
\label{fig:configuration}
\end{figure}

An example involving the simulation of a tri-variate stochastic vector process representing wind turbulent velocity fluctuations is provided here to illustrate the application of the proposed methodology. This example is modified from \cite{Deodatis1996}. Consider three components of the simulated vector process denoted by $f_{1}(t), f_{2}(t), f_{3}(t)$, describing the wind velocity fluctuations at three vertical points in a wind profile (points 1,2 and 3 in Figure \ref{fig:configuration}). 

The components of the $2^{nd}$-order cross spectrum (cross power spectrum) are given by
\begin{equation}
\begin{aligned}
	& S_{jj}(\omega) = S_{j}(\omega) \ j=1,2,3 \\
	& S_{jk}(\omega) = \sqrt{S_{j}(\omega)S_{k}(\omega)}\gamma_{jk}(\omega) \ j,k=1,2,3 \ j \neq k 
\end{aligned}
\end{equation}
where $S_{j}(\omega)$ is the power spectrum of process $f_{j}(t)$ and $\gamma_{jk}(\omega)$ is the coherence function between processes $f_{j}(t)$ and $f_{k}(t)$. The form suggested by Kaimal\cite{Kaimal1972} is selected to model the $2^{nd}$-order cross spectrum of the wind fluctuations and is given by
\begin{equation}
\begin{aligned}
	& S(z, \omega) = \frac{1}{2} \frac{200}{2\pi} u_{*}^{2} \frac{z}{U(z)} \frac{1}{[1 + 50 \frac{\omega z}{2 \pi U(z)}]^{\frac{5}{2}}} 
\end{aligned}
\end{equation}
where $z$ = height above the ground (in m); $u_{*}$ = shear velocity of the flow (in m/s); and $U(z)$ = mean wind speed at the height $z$ (in m/s). The model suggested in Davenport \cite{Davenport1968} is selected for the coherence function between the wind velocity fluctuations at different heights given by:
\begin{equation}
\begin{aligned}
	& \gamma(\Delta z, \omega) = \exp \Big[ \frac{-\omega}{2\pi} \frac{C_{z}\Delta z}{\frac{1}{2}[U(z_{1}) + U(z_{2})]} \Big] 
\end{aligned}
\end{equation}
where $U(z_{1})$ and $U(z_{2})$ are the mean wind speeds at heights $z_{1}$ and $z_{2}$ respectively, $\Delta z = |z_{1} - z_{2}|$, and $C_{z}$ is a constant equal to 10 for structural applications.

\begin{figure}[ht!]
\centering
  \includegraphics[width=0.8\linewidth]{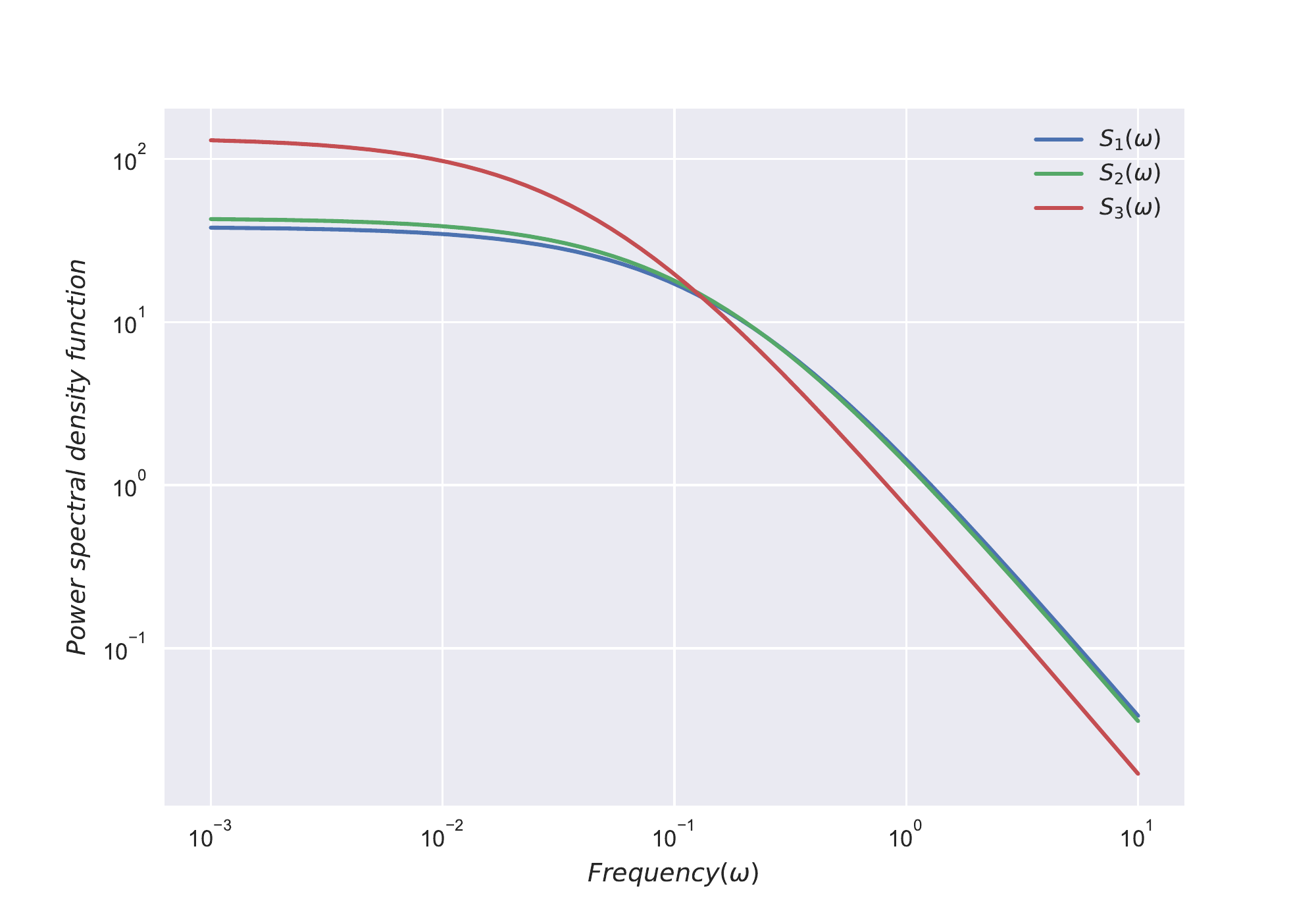}
  \caption{Power spectral density functions ($S_{j}(\omega); j=1,2,3$) for each component of the wind velocity vector.}
  \label{fig:all_s}
\end{figure}

\begin{figure}[ht!]
\centering
  \includegraphics[width=0.8\linewidth]{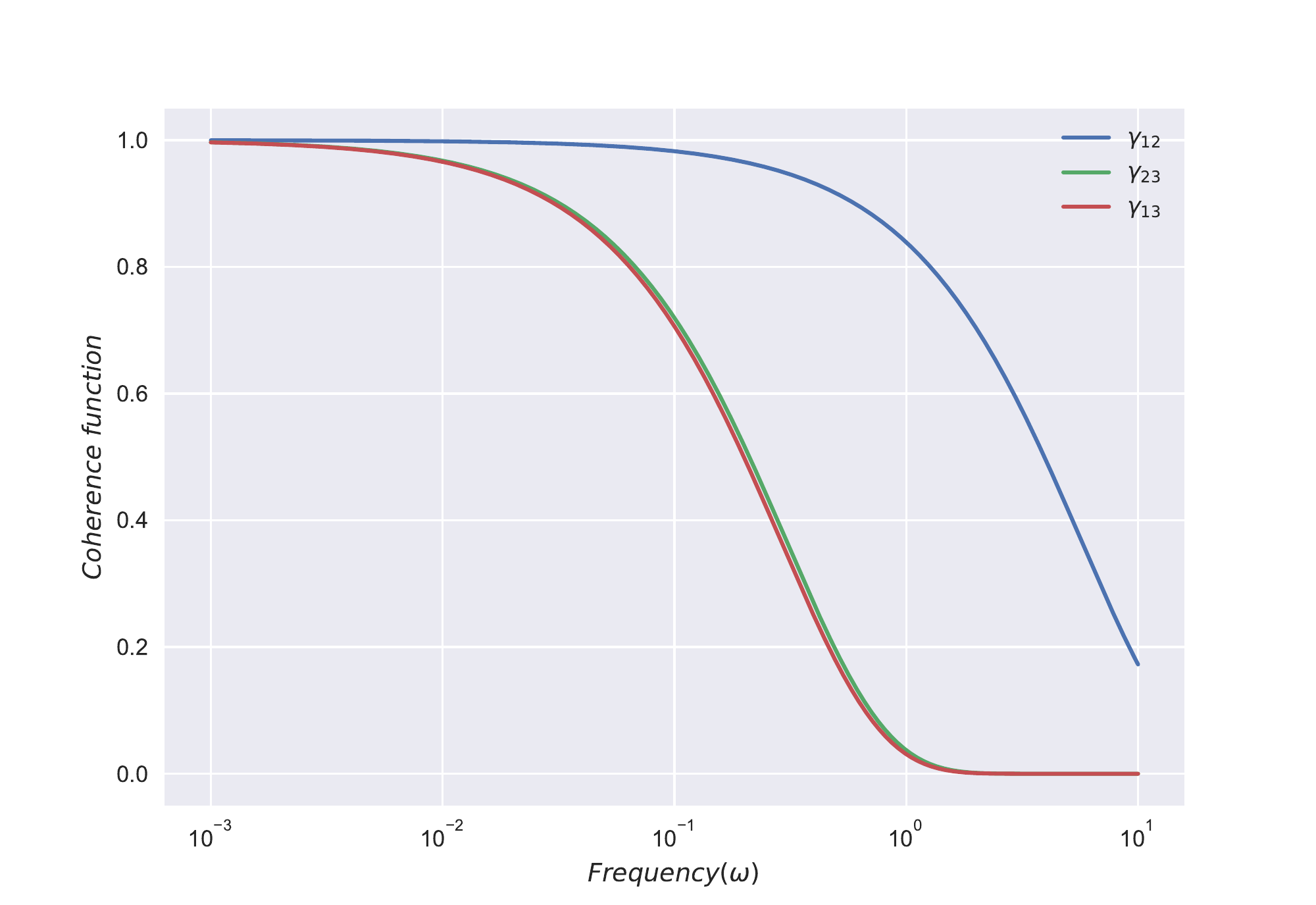}
  \caption{Coherence functions ($\gamma_{jk}(\omega); j,k=1,2,3, j \neq k$) between the wind velocity vector components.}
  \label{fig:all_gamma}
\end{figure}

The specific values of various parameters are obtained from \cite{Deodatis1996} and using the numerical results, the elements of the cross power spectral density are given by
\begin{equation}
\begin{aligned}
    &S_{11} = \frac{38.3}{(1 + 6.19 * w)^\frac{5}{3}}\\
    &S_{22} = \frac{43.3}{(1 + 6.98 * w)^\frac{5}{3}}\\
    &S_{33} = \frac{135}{(1 + 21.8 * w)^\frac{5}{3}}
\end{aligned}
\end{equation}
and the corresponding coherence functions are given as
\begin{equation}
\begin{aligned}
    &\gamma_{12}(\omega) = e^{-0.1757\omega}\\
    &\gamma_{13}(\omega) = e^{-3.478\omega}\\
    &\gamma_{23}(\omega) = e^{-3.292\omega}
\end{aligned}
\end{equation}
The spectra and coherences can be visualized in Fig. \ref{fig:all_s} and Fig. \ref{fig:all_gamma}.

The diagonal components of the 3rd-order cross-spectrum (cross-bispectrum) are assumed to take the following form:
\begin{equation}
\begin{aligned}
    &B_{111}(w_{1}, w_{2}) = \frac{50}{(1 + 6.19 * (w_{1} + w_{2}))^{\frac{5}{3}}}\\
    &B_{222}(w_{1}, w_{2}) = \frac{50}{(1 + 6.98 * (w_{1} + w_{2}))^{\frac{5}{3}}}\\
    &B_{333}(w_{1}, w_{2}) = \frac{50}{(1 + 21.8 * (w_{1} + w_{2}))^{\frac{5}{3}}}
\end{aligned}
\end{equation}
while the off-diagonal terms are given by
\begin{equation}
\begin{aligned}
    &B_{ijk}(w_{1}, w_{2}) = \sqrt[3]{B_{iii}(w_{1}, w_{2})B_{jjj}(w_{1}, w_{2})B_{kkk}(w_{1}, w_{2})}\gamma_{ijk}
\end{aligned}
\end{equation}
where $\gamma_{ijk}$ are the third-order coherence functions (or bi-coherences). The third-order coherence functions are given by
\begin{equation}
\begin{aligned}
    &\gamma_{112}(w_{1}, w_{2}) = e^{-0.171(w_{1} + w_{2})}\\
    &\gamma_{122}(w_{1}, w_{2}) = e^{-0.357(w_{1} + w_{2})}\\
    &\gamma_{113}(w_{1}, w_{2}) = e^{-1.287(w_{1} + w_{2})}\\
    &\gamma_{133}(w_{1}, w_{2}) = e^{-1.589(w_{1} + w_{2})}\\
    &\gamma_{123}(w_{1}, w_{2}) = e^{-3.473(w_{1} + w_{2})}\\
    &\gamma_{223}(w_{1}, w_{2}) = e^{-2.659(w_{1} + w_{2})}\\
    &\gamma_{233}(w_{1}, w_{2}) = e^{-2.775(w_{1} + w_{2})}
\end{aligned}
\end{equation}


Sample functions of this tri-variate stochastic wind velocity process are simulated using Eqs.\ \eqref{eqn:FFT1} and \eqref{eqn:FFT2} with the FFT technique. The upper cutoff frequency and the number of frequency discretizations are given by
\begin{equation}
    \omega_{u} = 2 \ \text{rad}/s; \ N_{\omega} = 100\\
\end{equation}
which results in the following frequency and time discretizations:
\begin{equation}
    \Delta \omega = 0.02 \ \text{rad}/s; \ \Delta t = 1.57 \ sec; \ T_{0} = 314.15 \ sec
\end{equation}

A single realization of each of the vector components is plotted in Fig. \ref{fig:var_1}, which also shows comparisons with sample functions generated using the 2nd-order SRM (having the same random phase angles) for comparison. 
\begin{figure}[!ht]
\centering
\begin{subfigure}
\centering
  \includegraphics[width=0.55\linewidth]{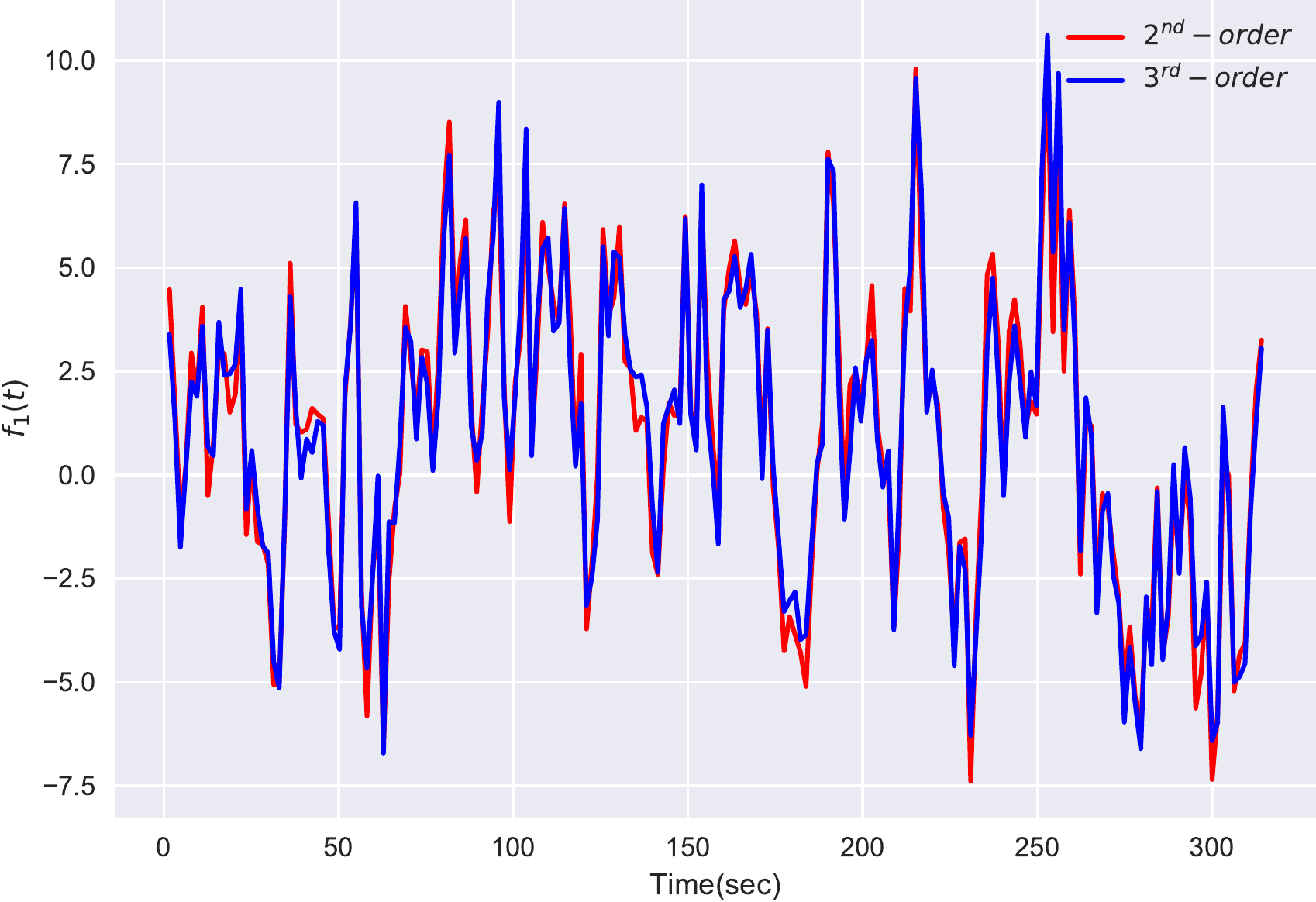}
  \label{fig:var_1}
\end{subfigure}
\begin{subfigure}
\centering
  \includegraphics[width=0.55\linewidth]{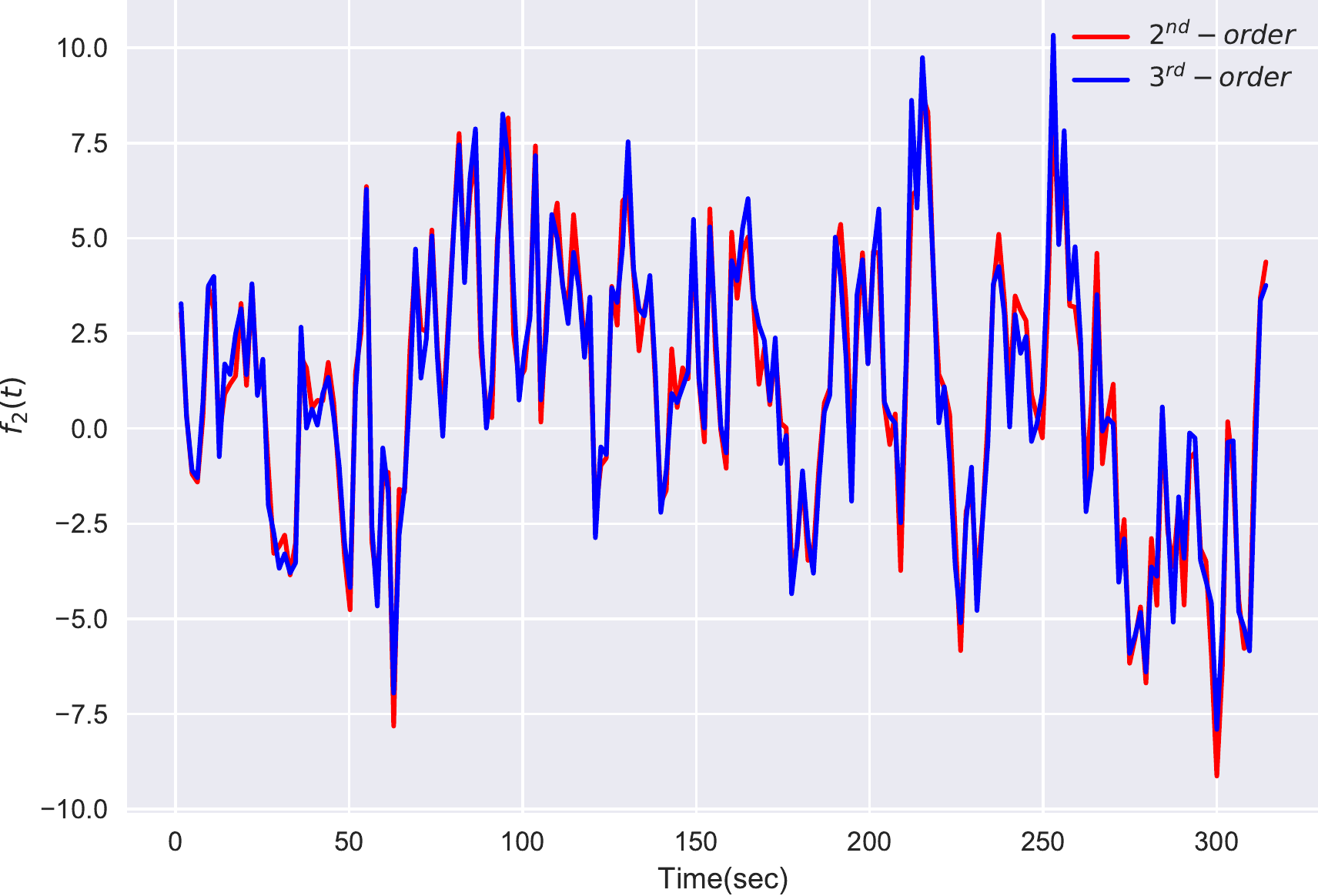}
  \label{fig:var_1}
\end{subfigure}
\begin{subfigure}
\centering
  \includegraphics[width=0.55\linewidth]{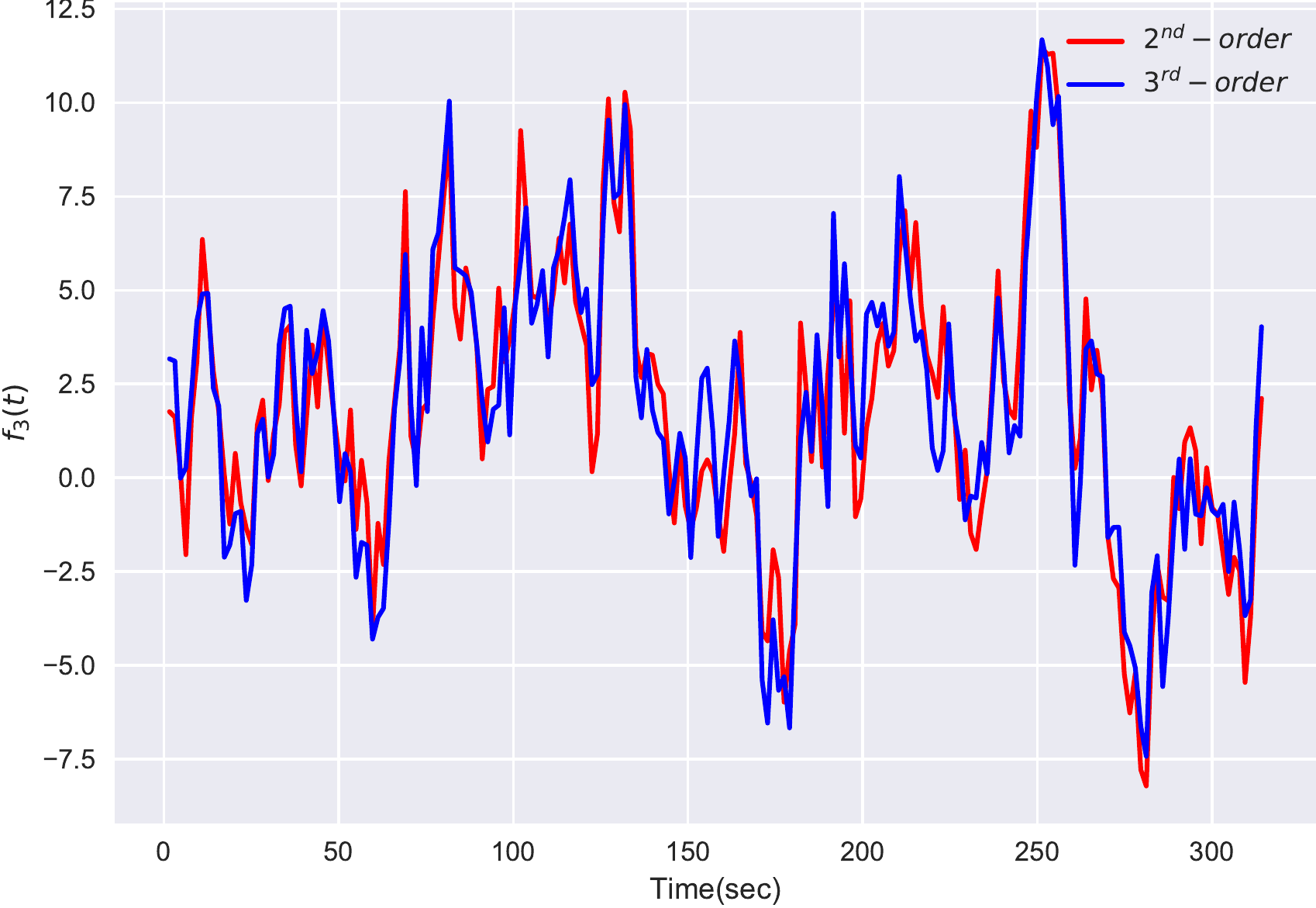}
  \label{fig:var_3}
\end{subfigure}
\caption{Velocity histories at points 1 (top), 2 (middle) and 3 (bottom).}
\end{figure}



To verify that the simulations are, indeed, possessing the prescribed statistical properties, the $1^{st}$, $2^{nd}$ and $3^{rd}$-order ensemble properties of the stochastic vector process are summarised in Tables \ref{table:example_1st_order} -- \ref{table:example_3rd_order}. Note that we do not produce plots of the spectral quantities because these are difficult to visualize.
\begin{table}[!ht]
\centering
\begin{tabular}{l l l l}
\hline
\textbf{moments} & \textbf{$3^{rd}$-order} & \textbf{$2^{nd}$-order} & \textbf{Target}\\
\hline
$\mathbb{E}[f_{1}(t)]$ & -0.00143 & -0.00143 & 0.00\\
$\mathbb{E}[f_{2}(t)]$ & -0.00147 & -0.00147 & 0.00\\
$\mathbb{E}[f_{3}(t)]$ & -0.00279 & -0.00279 & 0.00\\
\hline
\end{tabular}
\caption{First order statistics of the simulated vector process}
\label{table:example_1st_order}
\end{table}

\begin{table}[!ht]
\centering
\begin{tabular}{l l l l}
\hline
\textbf{moments} & \textbf{$3^{rd}$-order} & \textbf{$2^{nd}$-order} & \textbf{Target}\\
\hline
$\mathbb{E}[f_{1}^{2}(t)]$ & 14.541 & 14.538 & 14.539\\
$\mathbb{E}[f_{2}^{2}(t)]$ & 14.722 & 14.720 & 14.722\\
$\mathbb{E}[f_{3}^{2}(t)]$ & 14.724 & 14.723 & 14.723\\
$\mathbb{E}[f_{1}(t)f_{2}(t)]$ & 13.698 & 13.697 & 13.698\\
$\mathbb{E}[f_{1}(t)f_{3}(t)]$ & 7.628 & 7.627 & 7.628\\
$\mathbb{E}[f_{2}(t)f_{3}(t)]$ & 8.006 & 8.004 & 8.005\\
\hline
\end{tabular}
\caption{Second order statistics of the simulated vector process}
\label{table:example_2nd_order}
\end{table}

\begin{table}[!ht]
\centering
\begin{tabular}{l l l l}
\hline
\textbf{moments} & \textbf{$3^{rd}$-order} & \textbf{$2^{nd}$-order} & \textbf{Target}\\
\hline
$\mathbb{E}[f_{1}^{3}(t)]$ & 4.880 & 0.012 & 4.801\\
$\mathbb{E}[f_{2}^{3}(t)]$ & 3.870 & 0.004 & 3.825\\
$\mathbb{E}[f_{3}^{3}(t)]$ & 0.337 & -0.010 & 0.368\\
$\mathbb{E}[f_{1}(t)f_{2}^{2}(t)]$ & 3.938 & -0.051 & 3.939\\
$\mathbb{E}[f_{1}(t)f_{3}^{2}(t)]$ & 3.218 & -0.054 & 3.231\\
$\mathbb{E}[f_{2}(t)f_{1}^{2}(t)]$ & 0.931 & -0.048 & 0.981\\
$\mathbb{E}[f_{2}(t)f_{3}^{2}(t)]$ & 0.297 & -0.058 & 0.391\\
$\mathbb{E}[f_{3}(t)f_{1}^{2}(t)]$ & 0.435 & -0.057 & 0.513\\
$\mathbb{E}[f_{3}(t)f_{2}^{2}(t)]$ & 0.140 & -0.065 & 0.247\\
$\mathbb{E}[f_{1}(t)f_{2}(t)f_{3}(t)]$ & 0.355 & -0.050 & 0.425\\
\hline
\end{tabular}
\caption{Third order statistics of the simulated vector process}
\label{table:example_3rd_order}
\end{table}
From these tables, we see that the first and second-order ensemble moments are very close to the target for both the second and third-order simulations. However, the second-order simulations cannot match the target third-order moments (Table \ref{table:example_3rd_order}). The third-order simulations, on the other hand, match all moments up to third-order with very high accuracy.

\section{Conclusions}
\label{S:Conclusion}

In this, the second part of a two-part paper on the third-order spectral representation method (SRM), we present the simulation of stationary and ergodic, third-order multi-variate stochastic processes with the use of the fast Fourier transform. We first introduce a frequency multi-indexing for the uni-variate third-order SRM that yields ergodic sample functions. We prove the ergodicity of these simulated sample functions. Next, the simulation method is generalized for the case of multi-variate stochastic processes and proofs of ensemble and ergodic properties of the sample functions are provided in appendices. Finally, the simulation equation is reformulated in a way that allows the use of the fast Fourier transform to significantly improve computational efficiency. The method is then applied for the simulation of a tri-variate stochastic wind velocity field.

\section{Acknowledgements}
This work has been supported by the National Science Foundation under award number 1652044.

\bibliographystyle{model1-num-names}
\bibliography{elsarticle-template-1-num.bib}

\appendix

\section{Ensemble properties of simulated stochastic vector processes}
\label{A:1}

In this appendix, we prove that the sample functions generated by the third-order spectral representation method for stochastic vector processes match the target properties up to third-order in ensemble.

\subsection{First Order Properties}

The multi-variate stochastic vector processes have zero mean in ensemble. This proof is trivial and can be easily derived from commutative property of expectation.


\subsection{Second Order Properties}

The second-order cross-correlation function of the stochastic vector processes is shown to be equal to the prescribed second-order cross-correlation function as follows.

\begin{equation}
\begin{aligned}
    & R_{ab}(\tau) = \mathbb{E}[(f_{a}(t)f_{b}(t+\tau))] \\
    & = 4 \sum_{k_{1}=0}^{N-1} \sum_{k_{2}=0}^{N-1} \mathbb{E} \big[\sum_{l_{1}=1}^{m}\sum_{l_{2}=1}^{m}|H_{al_{1}}(\omega_{l_{1}k_{1}})||H_{bl_{2}}(\omega_{l_{2}k_{2}})| \Delta \omega\\
    & \cos(\omega_{k_{1}}t - \theta_{al_{1}}(\omega_{l_{1}k_{1}}) + \phi_{l_{1}k_{1}})\cos(\omega_{k_{2}}(t + \tau) - \theta_{al_{2}}(\omega_{l_{2}k_{2}}) + \phi_{l_{2}k_{2}}) \\
    & + \sum_{l_{1}=1}^{m}\sum_{l_{2}=1}^{m}\sum_{n_{2}=1}^{m}\sum_{p_{2}=1}^{m}\sum_{q_{2}=1}^{m}\sum_{i_{2} + j_{2} = k_{2}}^{i_{2} \geq j_{2}\geq 0} |H_{al_{1}}(\omega_{l_{1}k_{1}})| |B_{bl_{2}n_{2}}(\omega_{p_{2}i_{2}}, \omega_{q_{2}j_{2}})|\\
    & |G_{l_{2}p_{2}}(\omega_{p_{2}i_{2}})||G_{n_{2}q_{2}}(\omega_{q_{2}j_{2}})| \Delta \omega\cos(\omega_{l_{1}k_{1}}t - \theta_{al_{1}}(\omega_{l_{1}k_{1}}) + \phi_{l_{1}k_{1}}) \\
    & \cos((\omega_{p_{2}i_{2}} + \omega_{p_{2}j_{2}})(t+\tau) - \beta_{bl_{2}n_{2}}(\omega_{p_{2}i_{2}}, \omega_{q_{2}j_{2}}) - \theta^{I}_{l_{2}p_{2}}(\omega_{p_{2}i_{2}}) - \theta^{I}_{n_{2}q_{2}}(\omega_{q_{2}j_{2}}) + \phi_{p_{2}i_{2}} + \phi_{q_{2}j_{2}})\\
& + \sum_{l_{2}=1}^{m}\sum_{l_{1}=1}^{m}\sum_{n_{1}=1}^{m}\sum_{p_{1}=1}^{m}\sum_{q_{1}=1}^{m}\sum_{i_{1} + j_{1} = k_{1}}^{i_{1} \geq j_{1}\geq 0} |H_{al_{2}}(\omega_{l_{2}k_{2}})| |B_{al_{1}n_{1}}(\omega_{i_{1}}, \omega_{j_{1}})|\\
    & |G_{l_{1}p_{1}}(\omega_{p_{1}i_{1}})||G_{n_{1}q_{1}}(\omega_{q_{1}j_{1}})| \Delta \omega\cos(\omega_{l_{2}k_{2}}(t + \tau) - \theta_{al_{2}}(\omega_{l_{2}k_{2}}) + \phi_{l_{2}k_{2}}) \\
    & \cos((\omega_{p_{1}i_{1}} + \omega_{q_{1}j_{1}})t - \beta_{al_{1}n_{1}}(\omega_{p_{1}i_{1}}, \omega_{q_{1}j_{1}}) - \theta^{I}_{l_{1}p_{1}}(\omega_{p_{1}i_{1}}) - \theta^{I}_{n_{1}q_{1}}(\omega_{q_{1}j_{1}}) + \phi_{p_{1}i_{1}} + \phi_{q_{1}j_{1}})\\
    & + \sum_{l_{1}=1}^{m}\sum_{n_{1}=1}^{m}\sum_{p_{1}=1}^{m}\sum_{q_{1}=1}^{m}\sum_{l_{2}=1}^{m}\sum_{n_{2}=1}^{m}\sum_{p_{2}=1}^{m}\sum_{q_{2}=1}^{m}\sum_{i_{1} + j_{1} = k_{1}}^{i_{1} \geq j_{1}\geq 0}\sum_{i_{2} + j_{2} = k_{2}}^{i_{2} \geq j_{2}\geq 0}\\
    & |B_{al_{1}n_{1}}(\omega_{p_{1}i_{1}}, \omega_{q_{1}j_{1}})||G_{l_{1}p_{1}}(\omega_{p_{1}i_{1}})||G_{n_{1}q_{1}}(\omega_{q_{1}j_{1}})||B_{bl_{2}n_{2}}(\omega_{p_{2}i_{2}}, \omega_{q_{2}j_{2}})||G_{l_{2}p_{2}}(\omega_{p_{2}i_{2}})||G_{n_{2}q_{2}}(\omega_{q_{2}j_{2}})|\\
    & \cos((\omega_{p_{1}i_{1}} + \omega_{q_{1}j_{1}})t - \beta_{al_{1}n_{1}}(\omega_{p_{1}i_{1}}, \omega_{q_{1}j_{1}}) - \theta^{I}_{l_{1}p_{1}}(\omega_{p_{1}i_{1}}) - \theta^{I}_{n_{1}q_{1}}(\omega_{q_{1}j_{1}}) + \phi_{p_{1}i_{1}} + \phi_{q_{1}j_{1}})\\
    & \cos((\omega_{p_{2}i_{2}} + \omega_{q_{2}j_{2}})(t + \tau) - \beta_{bl_{2}n_{2}}(\omega_{p_{2}i_{2}}, \omega_{q_{2}j_{2}}) - \theta^{I}_{l_{2}p_{2}}(\omega_{p_{2}i_{2}}) - \theta^{I}_{n_{2}q_{2}}(\omega_{q_{2}j_{2}}) + \phi_{p_{2}i_{2}} + \phi_{q_{2}j_{2}}) \big] \\
\end{aligned}
\end{equation}

\noindent
The first term of this expression is given by
\begin{equation}
\begin{aligned}
    & 4\mathbb{E}\big[\sum_{k_{1}=0}^{N-1} \sum_{k_{2}=0}^{N-1}\sum_{l_{1}=1}^{m}\sum_{l_{2}=1}^{m}|H_{al_{1}}(\omega_{l_{1}k_{1}})||H_{bl_{2}}(\omega_{l_{2}k_{2}})| \Delta \omega\\
    & \cos(\omega_{l_{1}k_{1}}t - \theta_{al_{1}}(\omega_{l_{1}k_{1}}) + \phi_{l_{1}k_{1}})\cos(\omega_{l_{2}k_{2}}t - \theta_{al_{2}}(\omega_{l_{2}k_{2}}) + \phi_{l_{2}k_{2}}) \big] \\
\end{aligned}
\end{equation}
Since $\phi$ are random phase angles drawn from $[0, 2\pi]$, we have that when $k_{1} = k_{2} = k$ and $l_{1} = l_{2} = m$, this expression simplifies as
\begin{equation}
\begin{aligned}
    &4\mathbb{E}\big[\sum_{k=0}^{N-1}\sum_{m=1}^{m}|H_{am}(\omega_{mk})||H_{bm}(\omega_{mk})| \Delta \omega\cos(\omega_{mk}t - \theta_{am}(\omega_{mk}) + \phi_{mk}) \cos(\omega_{mk}(t + \tau) - \theta_{bm}(\omega_{mk}) + \phi_{mk}) \big] \\
    &= 2\mathbb{E}\big[\sum_{k=0}^{N-1}\sum_{m=1}^{m}|H_{am}(\omega_{mk})||H_{bm}(\omega_{mk})| \Delta \omega[\cos(\omega_{mk}\tau + \theta_{am}(\omega_{mk}) - \theta_{bm}(\omega_{mk}))\\
    & + \cos(\omega_{mk}(2t + \tau) - \theta_{am}(\omega_{mk}) - \theta_{bm}(\omega_{mk}) + 2\phi_{mk})] \big] \\
    &= 2\big[\sum_{k=0}^{N-1}\sum_{m=1}^{m}|H_{am}(\omega_{mk})||H_{bm}(\omega_{mk})| \Delta \omega\mathbb{E}[\cos(\omega_{mk}\tau + \theta_{am}(\omega_{mk}) - \theta_{bm}(\omega_{mk}))]\\
    & +\mathbb{E}[ \cos(\omega_{mk}(2t + \tau) - \theta_{am}(\omega_{mk}) - \theta_{bm}(\omega_{mk}) +  2\phi_{mk})] \big] \\
    &= 2\sum_{k=0}^{N-1}\sum_{m=1}^{m}|H_{am}(\omega_{mk})||H_{bm}(\omega_{mk})| \Delta \omega \cos(\omega_{mk}\tau + \theta_{am}(\omega_{mk}) - \theta_{bm}(\omega_{mk}))\\
    &= 2\int_{0}^{\omega_{u}}\sum_{m=1}^{m} |H_{am}(\omega_{m})||H_{bm}(\omega_{m})| e^{\iota(\omega_{m}\tau + \theta_{am}(\omega_{m}) - \theta_{b}(\omega_{m}))} d \omega\\
    &= \int_{-\omega_{u}}^{\omega_{u}} \sum_{m=1}^{m}H_{am}(\omega_{m})H^{*}_{mb}(\omega_{m}) e^{\iota\omega_{m}\tau} d \omega\\
    &= \int_{-\omega_{u}}^{\omega_{u}} S^{p}_{ab}(\omega) e^{\iota\omega\tau} d\omega\\
    &\approx \int_{-\infty}^{\infty} S^{p}_{ab}(\omega) e^{\iota\omega\tau} d\omega\\
\end{aligned}
\end{equation}

\noindent
The second and third terms of the expression are given by
\begin{equation}
\begin{aligned}
    & 4\sum_{k_{1}=0}^{N-1} \sum_{k_{2}=0}^{N-1} \mathbb{E} \big[ \sum_{l_{1}=1}^{m}\sum_{l_{2}=1}^{m}\sum_{n_{2}=1}^{m}\sum_{p_{2}=1}^{m}\sum_{q_{2}=1}^{m}\sum_{i_{2} + j_{2} = k_{2}}^{i_{2} \geq j_{2}\geq 0} |H_{al_{1}}(\omega_{l_{1}k_{1}})| |B_{bl_{2}n_{2}}(\omega_{p_{2}i_{2}}, \omega_{q_{2}j_{2}})|\\
    & |G_{l_{2}p_{2}}(\omega_{p_{2}i_{2}})||G_{n_{2}q_{2}}(\omega_{q_{2}j_{2}})| (\Delta\omega)^{\frac{3}{2}}\cos(\omega_{l_{1}k_{1}}t - \theta_{al_{1}}(\omega_{l_{1}k_{1}}) + \phi_{l_{1}k_{1}}) \\
    & \cos((\omega_{p_{2}i_{2}} + \omega_{q_{2}j_{2}})(t + \tau) - \beta_{al_{2}n_{2}}(\omega_{p_{2}i_{2}}, \omega_{q_{2}j_{2}}) - \theta^{I}_{al_{2}}(\omega_{p_{2}i_{2}}) - \theta^{I}_{an_{2}}(\omega_{q_{2}j_{2}}) + \phi_{p_{2}i_{2}} + \phi_{q_{2}j_{2}}) \big] \\
    & = 4 \sum_{k_{1}=0}^{N-1} \sum_{k_{2}=0}^{N-1} \mathbb{E} \big[\sum_{l_{2}=1}^{m}\sum_{l_{1}=1}^{m}\sum_{n_{1}=1}^{m}\sum_{p_{1}=1}^{m}\sum_{q_{1}=1}^{m}\sum_{i_{1} + j_{1} = k_{1}}^{i_{1} \geq j_{1}\geq 0} |H_{al_{2}}(\omega_{l_{2}k_{2}})| |B_{al_{1}n_{1}}(\omega_{p_{1}i_{1}}, \omega_{q_{1}j_{1}})|\\
    & |G_{l_{1}p_{1}}(\omega_{p_{1}i_{1}})||G_{n_{1}q_{1}}(\omega_{q_{1}j_{1}})|(\Delta\omega)^{\frac{3}{2}}\cos(\omega_{l_{2}k_{2}}(t + \tau) - \theta_{al_{2}}(\omega_{l_{2}k_{2}}) + \phi_{l_{2}k_{2}}) \\
    & \cos((\omega_{p_{1}i_{1}} + \omega_{q_{1}j_{1}})t - \beta_{al_{1}n_{1}}(\omega_{p_{1}i_{1}}, \omega_{q_{1}j_{1}}) - \theta^{I}_{l_{1}p_{1}}(\omega_{p_{1}i_{1}}) - \theta^{I}_{n_{1}q_{1}}(\omega_{q_{1}j_{1}}) + \phi_{p_{1}i_{1}} + \phi_{q_{1}j_{1}}) \big] = 0\\
\end{aligned}
\label{eqn:2_3_coupling_R_2}
\end{equation}

Both combinations of terms defined in \ref{eqn:2_3_coupling_R_2} are equal to zero in expectation since both involve coupling of odd number phase angles. 

Finally, the fourth term of the expression is given by
\begin{equation}
\begin{aligned}
    & 4\sum_{k_{1}=0}^{N-1} \sum_{k_{2}=0}^{N-1} \mathbb{E} \big[ \sum_{l_{1}=1}^{m}\sum_{n_{1}=1}^{m}\sum_{p_{1}=1}^{m}\sum_{q_{1}=1}^{m}\sum_{l_{2}=1}^{m}\sum_{n_{2}=1}^{m}\sum_{p_{2}=1}^{m}\sum_{q_{2}=1}^{m}\sum_{i_{1} + j_{1} = k_{1}}^{i_{1} \geq j_{1}\geq 0}\sum_{i_{2} + j_{2} = k_{2}}^{i_{2} \geq j_{2}\geq 0}\\
    & |B_{al_{1}n_{1}}(\omega_{p_{1}i_{1}}, \omega_{q_{1}j_{1}})| |G_{l_{1}p_{1}}(\omega_{p_{1}i_{1}})||G_{n_{1}q_{1}}(\omega_{q_{1}j_{1}})|\\
    & |B_{bl_{2}n_{2}}(\omega_{p_{2}i_{2}}, \omega_{q_{2}j_{2}})| |G_{l_{2}p_{2}}(\omega_{p_{2}i_{2}})||G_{n_{2}q_{2}}(\omega_{q_{1}j_{2}})|\Delta \omega^{2} \\
    & \cos((\omega_{p_{1}i_{1}} + \omega_{q_{1}j_{1}})t - \beta_{al_{1}n_{1}}(\omega_{p_{1}i_{1}}, \omega_{q_{1}j_{1}}) - \theta^{I}_{l_{1}p_{1}}(\omega_{p_{1}i_{1}}) - \theta^{I}_{n_{1}q_{1}}(\omega_{q_{1}j_{1}}) + \phi_{p_{1}i_{1}} + \phi_{q_{1}j_{1}})\\
    & \cos((\omega_{p_{2}i_{2}} + \omega_{q_{2}j_{2}})(t + \tau) - \beta_{bl_{2}n_{2}}(\omega_{p_{2}i_{2}}, \omega_{q_{2}j_{2}}) - \theta^{I}_{l_{2}p_{2}}(\omega_{p_{2}i_{2}}) - \theta^{I}_{n_{2}q_{2}}(\omega_{q_{2}j_{2}}) + \phi_{p_{2}i_{2}} + \phi_{q_{2}j_{2}}) \big]\\
\end{aligned}
\end{equation}
When $p_{1} = p_{2} = p$, $q_{1} = q_{2} = q$ and $i_{1} = i_{2} = i$, $j_{1} = j_{2} = j$, implying $k_{1} = k_{2} = k$, this simplifies as
\begin{equation}
\begin{aligned}
    &4\sum_{k=0}^{N-1}\mathbb{E}\big[ \sum_{l_{1}=1}^{m} \sum_{n_{1}=1}^{m} \sum_{l_{2}=1}^{m} \sum_{n_{2}=1}^{m} \sum_{p=1}^{m}\sum_{q=1}^{m}\sum_{i + j = k}^{i \geq j \geq 0} |B_{al_{1}n_{1}}(\omega_{pi}, \omega_{qj})||G_{l_{1}p}(\omega_{pi})||G_{n_{1}q}(\omega_{qj})|\\
    & |B_{bl_{2}n_{2}}(\omega_{pi}, \omega_{qj})||G_{l_{2}p}(\omega_{pi})||G_{n_{2}q}(\omega_{qj})|\Delta \omega^{2} \\
    & \cos((\omega_{pi} + \omega_{qj})t - \beta_{al_{1}n_{1}}(\omega_{pi}, \omega_{qj}) - \theta^{I}_{l_{1}p}(\omega_{pi}) - \theta^{I}_{n_{1}q}(\omega_{qj}) + \phi_{pi} + \phi_{qj})\\
    & \cos((\omega_{pi} + \omega_{qj})(t + \tau) - \beta_{bl_{2}n_{2}}(\omega_{pi}, \omega_{qj}) - \theta^{I}_{l_{2}p}(\omega_{pi}) - \theta^{I}_{n_{2}q}(\omega_{qj}) + \phi_{pi} + \phi_{qj})\big]\\
    &= 2\sum_{k=0}^{N-1}\mathbb{E}\big[ \sum_{l_{1}=1}^{m} \sum_{n_{1}=1}^{m} \sum_{l_{2}=1}^{m} \sum_{n_{2}=1}^{m} \sum_{p=1}^{m}\sum_{q=1}^{m}\sum_{i + j = k}^{i \geq j \geq 0} |B_{al_{1}n_{1}}(\omega_{pi}, \omega_{qj})||G_{l_{1}p}(\omega_{pi})||G_{n_{1}q}(\omega_{qj})|\\
    & |B_{bl_{2}n_{2}}(\omega_{pi}, \omega_{qj})||G_{l_{2}p}(\omega_{pi})||G_{n_{2}q}(\omega_{qj})|\Delta \omega^{2} \\
    & \cos(-(\omega_{pi} + \omega_{qj})t + \beta_{al_{1}n_{1}}(\omega_{pi}, \omega_{qj}) + \theta^{I}_{l_{1}p}(\omega_{pi}) + \theta^{I}_{n_{1}q}(\omega_{qj}) - \phi_{pi} - \phi_{qj}\\
    &  + (\omega_{pi} + \omega_{qj})(t + \tau) - \beta_{bl_{2}n_{2}}(\omega_{pi}, \omega_{qj}) - \theta^{I}_{l_{2}p}(\omega_{pi}) - \theta^{I}_{n_{2}q}(\omega_{qj}) + \phi_{pi} + \phi_{qj})\big]\\
    &= 2\sum_{k=0}^{N-1}\sum_{l_{1}=1}^{m}\sum_{n_{1}=1}^{m} \sum_{l_{2}=1}^{m}\sum_{n_{2}=1}^{m}\sum_{p=1}^{m}\sum_{q=1}^{m}\sum_{i + j = k}^{i \geq j \geq 0} |B_{al_{1}n_{1}}(\omega_{pi}, \omega_{qj})||G_{l_{1}p}(\omega_{pi})||G_{n_{1}q}(\omega_{qj})|\\
    & |B_{bl_{2}n_{2}}(\omega_{pi}, \omega_{qj})||G_{l_{2}p}(\omega_{pi})||G_{n_{2}q}(\omega_{qj})|\Delta \omega^{2} \\
    & \cos((\omega_{pi} + \omega_{qj})\tau + \beta_{al_{1}n_{1}}(\omega_{pi}, \omega_{qj}) + \theta^{I}_{l_{1}p}(\omega_{pi}) + \theta^{I}_{n_{1}q}(\omega_{qj})- \beta_{bl_{2}n_{2}}(\omega_{pi}, \omega_{qj}) - \theta^{I}_{l_{2}p}(\omega_{pi}) - \theta^{I}_{n_{2}q}(\omega_{qj}))\\
\end{aligned}
\end{equation}

\begin{equation*}
\begin{aligned}
    &= 2\sum_{k=0}^{N-1}\sum_{l_{1}=1}^{m}\sum_{n_{1}=1}^{m} \sum_{l_{2}=1}^{m}\sum_{n_{2}=1}^{m}\sum_{p=1}^{m}\sum_{q=1}^{m}\sum_{i + j = k}^{i \geq j \geq 0} |B_{al_{1}n_{1}}(\omega_{pi}, \omega_{qj})||G_{l_{1}p}(\omega_{pi})||G_{n_{1}q}(\omega_{qj})|\\
    & |B_{bl_{2}n_{2}}(\omega_{pi}, \omega_{qj})||G_{l_{2}p}(\omega_{pi})||G_{n_{2}q}(\omega_{qj})|\\
    & \cos((\omega_{pi} + \omega_{qj})\tau + \beta_{al_{1}n_{1}}(\omega_{pi}, \omega_{j}) + \theta^{I}_{l_{1}p}(\omega_{pi}) + \theta^{I}_{n_{1}q}(\omega_{qj})- \beta_{bl_{2}n_{2}}(\omega_{pi}, \omega_{j}) - \theta^{I}_{l_{2}p}(\omega_{pi}) - \theta^{I}_{n_{2}q}(\omega_{j}))\Delta\omega^{2}\\
    &= 2\sum_{k=0}^{N-1}\sum_{l_{1}=1}^{m}\sum_{n_{1}=1}^{m}\sum_{l_{2}=1}^{m}\sum_{n_{2}=1}^{m}\sum_{p=1}^{m}\sum_{q=1}^{m}\sum_{i + j = k}^{i \geq j \geq 0} |B_{al_{1}n_{1}}(\omega_{pi}, \omega_{qj})||G_{l_{1}p}(\omega_{pi})||G_{n_{1}q}(\omega_{qj})|\\
    & |B_{bl_{2}n_{2}}(\omega_{pi}, \omega_{qj})||G_{l_{2}p}(\omega_{pi})||G_{n_{2}q}(\omega_{qj})|e^{\iota(\omega_{pi} + \omega_{qj})\tau}e^{\iota\beta_{al_{1}n_{1}}(\omega_{pi}, \omega_{qj})}e^{\iota\theta^{I}_{l_{1}p}(\omega_{pi})}e^{\iota\theta^{I}_{n_{1}q}(\omega_{qj})}\\
    & e^{-\iota\beta_{bl_{2}n_{2}}(\omega_{pi}, \omega_{qj})}e^{-\iota\theta^{I}_{l_{2}p}(\omega_{pi})}e^{-\iota\theta^{I}_{n_{2}q}(\omega_{qj})}\Delta \omega^{2}\\
    &= 2\sum_{k=0}^{N-1}\sum_{l_{1}=1}^{m}\sum_{n_{1}=1}^{m}\sum_{l_{2}=1}^{m}\sum_{n_{2}=1}^{m}\sum_{p=1}^{m}\sum_{q=1}^{m}\sum_{i + j = k}^{i \geq j \geq 0} B_{al_{1}n_{1}}(\omega_{pi}, \omega_{qj})G_{l_{1}p}(\omega_{pi})G_{n_{1}q}(\omega_{qj})\\
    & B_{bl_{2}n_{2}}^{*}(\omega_{pi}, \omega_{qj})G_{l_{2}p}^{*}(\omega_{pi})G_{n_{2}q}^{*}(\omega_{qj})e^{\iota(\omega_{pi} + \omega_{qj})\tau}\Delta\omega^{2}\\
    &= 2\sum_{k=0}^{N-1}\sum_{l_{1}=1}^{m}\sum_{n_{1}=1}^{m}\sum_{l_{2}=1}^{m}\sum_{n_{2}=1}^{m}\sum_{p=1}^{m}\sum_{q=1}^{m}\sum_{i + j = k}^{i \geq j \geq 0} \\
    & B_{al_{1}n_{1}}(\omega_{pi}, \omega_{qj})B_{bl_{2}n_{2}}^{*}(\omega_{pi}, \omega_{qj})G_{l_{1}p}(\omega_{pi})G_{l_{2}p}^{*}(\omega_{pi})G_{n_{1}q}(\omega_{qj})G_{n_{2}q}^{*}(\omega_{qj})e^{\iota(\omega_{pi} + \omega_{qj})\tau}\Delta\omega^{2}\\
    & \approx 2\sum_{k=0}^{N-1}\sum_{i + j = k}^{i \geq j \geq 0}\sum_{l_{1}=1}^{m}\sum_{n_{1}=1}^{m}\sum_{l_{2}=1}^{m}\sum_{n_{2}=1}^{m} B_{al_{1}n_{1}}(\omega_{i}, \omega_{j})B_{bl_{2}n_{2}}^{*}(\omega_{i}, \omega_{j})S_{l_{1}l_{2}}^{(p)I}(\omega_{i})S_{n_{1}n_{2}}^{(p)I}(\omega_{i})e^{\iota\omega_{k}\tau}\Delta\omega^{2}\\
    & \approx 2\sum_{k=0}^{N-1} \Big( S_{ab}(\omega) - S_{ab}^{(p)}(\omega) \Big) e^{\iota\omega\tau}\Delta\omega\\
    & \approx 2\int_{0}^{\omega_{u}} \Big( S_{ab}(\omega) - S_{ab}^{(p)}(\omega) \Big) e^{\iota\omega\tau} d\omega\\
    & \approx \int_{-\omega_{u}}^{\omega_{u}} \Big( S_{ab}(\omega) - S_{ab}^{(p)}(\omega) \Big) e^{\iota\omega\tau} d\omega\\
    & \approx \int_{-\infty}^{\infty} \Big( S_{ab}(\omega) - S^{(p)}_{ab}(\omega)\Big) e^{\iota\omega\tau}d\omega\\
\end{aligned}
\end{equation*}

\noindent
Combining these four terms, we have
\begin{equation}
\begin{aligned}
    R_{ab}(\tau) &= \int_{-\infty}^{\infty} S^{(p)}_{ab}(\omega) e^{\iota\omega\tau}d\omega + 0 + 0 + \int_{-\infty}^{\infty} \Big(S_{ab}(\omega_{k}) - S^{(p)}_{ab}(\omega_{k})\Big) e^{\iota\omega_{k}\tau}d\omega\\
    R_{ab}(\tau) &= \int_{-\infty}^{\infty} S_{ab}(\omega) e^{\iota\omega\tau}d\omega\\
\end{aligned}
\end{equation}

\subsection{Third Order Properties}

The third-order cross-correlation function of the stochastic vector processes is shown to be equal to the prescribed third-order cross-correlation function as follows.
\begin{equation}
\begin{aligned}
    & R_{abc}(\tau_{1}, \tau_{2}) = \mathbb{E}[f_{a}(t)f_{b}(t+\tau_{1})f_{c}(t+\tau_{2})]\\
    & = 8 \mathbb{E} \big[ [\sum_{k_{1}=0}^{N-1} \big[\sum_{l_{1}=1}^{m}|H_{al_{1}}(\omega_{k_{1}})|\sqrt{\Delta \omega}\cos(\omega_{k_{1}}t - \theta_{al_{1}}(\omega_{k_{1}}) + \phi_{l_{1}k_{1}}) \\
    & + \sum_{l_{1}=1}^{m}\sum_{n_{1}=1}^{m}\sum_{p_{1}=1}^{m}\sum_{q_{1}=1}^{m}\sum_{i_{1} + j_{1} = k_{1}}^{i_{1} \geq j_{1}\geq 0} |B_{al_{1}n_{1}}(\omega_{i_{1}}, \omega_{j_{1}})||G_{l_{1}p_{1}}(\omega_{i_{1}})||G_{n_{1}q_{1}}(\omega_{j_{1}})|\Delta\omega\\
    & \cos((\omega_{i_{1}} + \omega_{j_{1}})t - \beta_{al_{1}n_{1}}(\omega_{i_{1}}, \omega_{j_{1}}) - \theta^{I}_{l_{1}p_{1}}(\omega_{i_{1}}) - \theta^{I}_{n_{1}q_{1}}(\omega_{j_{1}}) + \phi_{p_{1}i_{1}} + \phi_{q_{1}j_{1}})] \\
	& [\sum_{k_{2}=0}^{N-1} \big[\sum_{l_{2}=1}^{m}|H_{al_{2}}(\omega_{k_{2}})|\sqrt{\Delta \omega}\cos(\omega_{k_{2}}t - \theta_{al_{2}}(\omega_{k_{2}}) + \phi_{l_{2}k_{2}}) \\
    & + \sum_{l_{2}=1}^{m}\sum_{n_{2}=1}^{m}\sum_{p_{2}=1}^{m}\sum_{q_{2}=1}^{m}\sum_{i_{2} + j_{2} = k_{2}}^{i_{2} \geq j_{2}\geq 0} |B_{al_{2}n_{2}}(\omega_{i_{2}}, \omega_{j_{2}})||G_{l_{2}p_{2}}(\omega_{i_{2}})||G_{n_{2}q_{2}}(\omega_{j_{2}})|\Delta\omega\\
    & \cos((\omega_{i_{2}} + \omega_{j_{2}})t - \beta_{al_{2}n_{2}}(\omega_{i_{2}}, \omega_{j_{2}}) - \theta^{I}_{l_{2}p_{2}}(\omega_{i_{2}}) - \theta^{I}_{n_{2}q_{2}}(\omega_{j_{2}}) + \phi_{p_{2}i_{2}} + \phi_{q_{2}j_{2}})] \\
	& [\sum_{k_{3}=0}^{N-1} \big[\sum_{l_{3}=1}^{m}|H_{al_{3}}(\omega_{k_{3}})|\sqrt{\Delta \omega}\cos(\omega_{k_{3}}t - \theta_{al_{3}}(\omega_{k_{3}}) + \phi_{l_{3}k_{3}}) \\
    & + \sum_{l_{3}=1}^{m}\sum_{n_{3}=1}^{m}\sum_{p_{3}=1}^{m}\sum_{q_{3}=1}^{m}\sum_{i_{3} + j_{3} = k_{3}}^{i_{3} \geq j_{3}\geq 0} |B_{al_{3}n_{3}}(\omega_{i_{3}}, \omega_{j_{3}})||G_{l_{3}p_{3}}(\omega_{i_{3}})||G_{n_{3}q_{3}}(\omega_{j_{3}})|\Delta\omega\\
    & \cos((\omega_{i_{3}} + \omega_{j_{3}})t - \beta_{al_{3}n_{3}}(\omega_{i_{3}}, \omega_{j_{3}}) - \theta^{I}_{l_{3}p_{3}}(\omega_{i_{3}}) - \theta^{I}_{n_{3}q_{3}}(\omega_{j_{3}}) + \phi_{p_{3}i_{3}} + \phi_{q_{3}j_{3}})] \big] \\
\end{aligned}
\end{equation}
When this product is expanded, there are 6 symmetric components each taking the form
\begin{equation}
\begin{aligned}
	& 8\sum_{k_{1} = 1}^{N}\sum_{i_{1} + j_{1}= k_{1}}^{i_{1} \geq j_{1} \geq 0}\sum_{k_{2} = 1}^{N}\sum_{k_{3} = 1}^{N}\sum_{l_{1} = 1}^{m}\sum_{n_{1} = 1}^{m}\sum_{p_{1} = 1}^{m}\sum_{q_{1} = 1}^{m}\sum_{l_{2} = 1}^{m}\sum_{l_{3} = 1}^{m} \\
	& \mathbb{E} \big[ |B_{al_{1}n_{1}}(\omega_{i_{1}}, \omega_{j_{1}})||G_{l_{1}p_{1}}(\omega_{i_{1}})||G_{n_{1}q_{1}}(\omega_{j_{1}})| \Delta \omega \\
    & \cos((\omega_{i_{1}} + \omega_{j_{1}})t - \beta_{al_{1}n_{1}}(\omega_{i_{1}}, \omega_{j_{1}}) - \theta^{p}_{l_{1}p_{1}}(\omega_{i_{1}}) - \theta^{p}_{n_{1}q_{1}}(\omega_{j_{1}}) + \phi_{p_{1}i_{1}} + \phi_{q_{1}j_{1}}) \\
    & |H_{bl_{2}}(\omega_{k_{2}})|\sqrt{\Delta \omega}\cos(\omega_{k_{2}}(t + \tau_{1}) - \theta^{p}_{bl_{2}}(\omega_{k_{2}}) + \phi_{l_{2}k_{2}})\\
    & |H_{cl_{3}}(\omega_{k_{3}})|\sqrt{\Delta \omega}\cos(\omega_{k_{3}}(t + \tau_{2}) - \theta^{p}_{cl_{3}}(\omega_{k_{3}}) + \phi_{l_{3}k_{3}}) \big]\\
\end{aligned}
\end{equation}
Therefore we can express $R_{abc}(\tau_{1}, \tau_{2})$ as
\begin{equation}
\begin{aligned}
	& 48\sum_{k_{1} = 1}^{N}\sum_{i_{1} + j_{1}= k_{1}}^{i_{1} \geq j_{1} \geq 0}\sum_{k_{2} = 1}^{N}\sum_{k_{3} = 1}^{N}\sum_{l_{1} = 1}^{m}\sum_{n_{1} = 1}^{m}\sum_{p_{1} = 1}^{m}\sum_{q_{1} = 1}^{m}\sum_{l_{2} = 1}^{m}\sum_{l_{3} = 1}^{m} \mathbb{E} \big[ |B_{al_{1}n_{1}}(\omega_{i_{1}}, \omega_{j_{1}})||G_{l_{1}p_{1}}(\omega_{i_{1}})||G_{n_{1}q_{1}}(\omega_{j_{1}})| \Delta \omega \\
    & \cos((\omega_{i_{1}} + \omega_{j_{1}})t - \beta_{al_{1}n_{1}}(\omega_{i_{1}}, \omega_{j_{1}}) - \theta^{p}_{l_{1}p_{1}}(\omega_{i_{1}}) - \theta^{p}_{n_{1}q_{1}}(\omega_{j_{1}}) + \phi_{p_{1}i_{1}} + \phi_{q_{1}j_{1}}) \\
    & |H_{bl_{2}}(\omega_{k_{2}})|\sqrt{\Delta \omega}\cos(\omega_{k_{2}}(t + \tau_{1}) - \theta^{p}_{bl_{2}}(\omega_{k_{2}}) + \phi_{l_{2}k_{2}})\\
    & |H_{cl_{3}}(\omega_{k_{3}})|\sqrt{\Delta \omega}\cos(\omega_{k_{3}}(t + \tau_{2}) - \theta^{p}_{cl_{3}}(\omega_{k_{3}}) + \phi_{l_{3}k_{3}}) \big] \\
	& = 48\sum_{k_{1} = 1}^{N}\sum_{i_{1} + j_{1}= k_{1}}^{i_{1} \geq j_{1} \geq 0}\sum_{k_{2} = 1}^{N}\sum_{k_{3} = 1}^{N}\sum_{l_{1} = 1}^{m}\sum_{n_{1} = 1}^{m}\sum_{p_{1} = 1}^{m}\sum_{q_{1} = 1}^{m}\sum_{l_{2} = 1}^{m}\sum_{l_{3} = 1}^{m} \\
	& \mathbb{E}\big[|B_{al_{1}n_{1}}(\omega_{i_{1}}, \omega_{j_{1}})||G_{l_{1}p_{1}}(\omega_{i_{1}})||H_{bl_{2}}(\omega_{k_{2}})||G_{n_{1}q_{1}}(\omega_{j_{1}})||H_{cl_{3}}(\omega_{k_{3}})|\Delta \omega^{2} \\
    & \cos((\omega_{i_{1}} + \omega_{j_{1}})t - \beta_{al_{1}n_{1}}(\omega_{i_{1}}, \omega_{j_{1}}) - \theta^{p}_{l_{1}p_{1}}(\omega_{i_{1}}) - \theta^{p}_{n_{1}q_{1}}(\omega_{j_{1}}) + \phi_{p_{1}i_{1}} + \phi_{q_{1}j_{1}}) \\
    & \cos(\omega_{k_{2}}(t + \tau_{1}) - \theta^{p}_{bl_{2}}(\omega_{k_{2}}) + \phi_{l_{2}k_{2}})\cos(\omega_{k_{3}}(t + \tau_{2}) - \theta^{p}_{cl_{3}}(\omega_{k_{3}}) + \phi_{l_{3}k_{3}}) \big] \\
\end{aligned}
\end{equation}

\noindent
When $i_{1} = k_{2} = i$, $j_{1} = k_{3} = j$, implying $k_{2} + k_{3} = k_{1} = k$; $p_{1} = l_{2} = p$ and $q_{1} = l_{3} = q$, this expression simplifies as
\begin{equation}
\begin{aligned}
	& 48\sum_{k = 1}^{N}\sum_{i + j= k}^{i \geq j \geq 0}\sum_{m = 1}^{m}\sum_{n = 1}^{m}\sum_{p = 1}^{m}\sum_{q = 1}^{m}\mathbb{E}\big[|B_{amn}(\omega_{i}, \omega_{j})||G_{mp}(\omega_{i})||H_{bp}(\omega_{i})||G^{(p)}_{nq}(\omega_{j})||H_{cq}(\omega_{j})|\Delta \omega^{2} \\
    & \cos((\omega_{i} + \omega_{j})t - \beta_{amn}(\omega_{i}, \omega_{j}) - \theta^{p}_{mp}(\omega_{i}) - \theta^{p}_{nq}(\omega_{j}) + \phi_{pi} + \phi_{qj}) \\
    & \cos(\omega_{i}(t + \tau_{1}) - \theta^{p}_{bp}(\omega_{i}) + \phi_{pi})\cos(\omega_{j}(t + \tau_{2}) - \theta^{p}_{cq}(\omega_{j}) + \phi_{qj}) \big] \\
	& = 12\sum_{k = 1}^{N}\sum_{i + j= k}^{i \geq j \geq 0}\sum_{m = 1}^{m}\sum_{n = 1}^{m}\sum_{p = 1}^{m}\sum_{q = 1}^{m}\mathbb{E}\big[|B_{amn}(\omega_{i}, \omega_{j})||G_{mp}(\omega_{i})||H_{bp}(\omega_{i})||G_{nq}(\omega_{j})||H_{cq}(\omega_{j})|\Delta \omega^{2} \\
    & \cos(-(\omega_{i} + \omega_{j})t + \beta_{amn}(\omega_{i}, \omega_{j}) + \theta^{p}_{mp}(\omega_{i}) + \theta^{p}_{nq}(\omega_{j}) - \phi_{pi} - \phi_{qj} \\
    & + \omega_{i}(t + \tau_{1}) - \theta^{p}_{bp}(\omega_{i}) + \phi_{pi} + \omega_{j}(t + \tau_{2}) - \theta^{p}_{cq}(\omega_{j}) + \phi_{qj}) \big] \\
	& = 12\sum_{k = 1}^{N}\sum_{i + j= k}^{i \geq j \geq 0}\sum_{m = 1}^{m}\sum_{n = 1}^{m}\sum_{p = 1}^{m}\sum_{q = 1}^{m}|B_{amn}(\omega_{i}, \omega_{j})||G_{mp}(\omega_{i})||H_{bp}(\omega_{i})||G_{nq}(\omega_{j})||H_{cq}(\omega_{j})|\Delta \omega^{2} \\
    & \cos(\beta_{amn}(\omega_{i}, \omega_{j}) + \theta^{p}_{mp}(\omega_{i}) + \theta^{p}_{nq}(\omega_{j}) + \omega_{i}\tau_{1} - \theta^{p}_{bp}(\omega_{i}) + \omega_{j}\tau_{2} - \theta^{p}_{cq}(\omega_{j}))\\
\end{aligned}
\end{equation}

\begin{equation*}
\begin{aligned}
	& = 12\sum_{k = 1}^{N}\sum_{i + j= k}^{i \geq j \geq 0}\sum_{m = 1}^{m}\sum_{n = 1}^{m}\sum_{p = 1}^{m}\sum_{q = 1}^{m}|B_{amn}(\omega_{i}, \omega_{j})||G_{mp}(\omega_{i})||H_{bp}(\omega_{i})||G_{nq}(\omega_{j})||H_{cq}(\omega_{j})|\Delta \omega^{2} \\
    & e^{\iota\beta_{amn}(\omega_{i}, \omega_{j})} e^{\iota\theta^{p}_{mp}(\omega_{i})} e^{\iota \theta^{p}_{nq}(\omega_{j})} e^{\iota\omega_{i}\tau_{1}} e^{-\iota\theta^{p}_{bp}(\omega_{i})} e^{\iota\omega_{j}\tau_{2}} e^{-\iota\theta^{p}_{cq}(\omega_{j})}\\
	& = 12\sum_{k = 1}^{N}\sum_{i + j= k}^{i \geq j \geq 0}\sum_{m = 1}^{m}\sum_{n = 1}^{m}\sum_{p = 1}^{m}\sum_{q = 1}^{m}B_{amn}(\omega_{i}, \omega_{j})G_{mp}(\omega_{i})H^{(p)*}_{bp}(\omega_{i})G_{nq}(\omega_{j})H_{cq}(\omega_{j})\Delta \omega^{2} e^{\iota\omega_{i}\tau_{1}} e^{\iota\omega_{j}\tau_{2}}\\
	& = 12\sum_{k = 1}^{N}\sum_{i + j= k}^{i \geq j \geq 0}\sum_{m = 1}^{m}\sum_{n = 1}^{m}B_{al_{1}n_{1}}(\omega_{i}, \omega_{j})I_{bl_{1}}(\omega_{i}) I_{cn_{1}}(\omega_{j})\Delta \omega^{2} e^{\iota\omega_{i}\tau_{1}} e^{\iota\omega_{j}\tau_{2}}\\
	& = 12\sum_{k = 1}^{N}\sum_{i + j= k}^{i \geq j \geq 0} B_{abc}(\omega_{i}, \omega_{j})e^{\iota\omega_{i}\tau_{1}} e^{\iota\omega_{j}\tau_{2}}\Delta \omega^{2}\\
	& = 6 \int_{0}^{\infty} \int_{0}^{\infty} B_{abc}(\omega_{i}, \omega_{j}) e^{\iota\omega_{i}\tau_{1}} e^{\iota\omega_{j}\tau_{2}} d\omega_{i}d\omega_{j}\\
\end{aligned}
\end{equation*}
Hence,
\begin{equation}
\begin{aligned}
   R_{abc}(\tau_{1}, \tau_{2}) = \int_{-\infty}^{\infty} \int_{-\infty}^{\infty}B_{abc}(\omega_{i}, \omega_{j})e^{\iota\omega_{i}\tau_{1}}e^{\iota\omega_{j}\tau_{2}}d\omega_{i}d\omega_{j}
\end{aligned}
\end{equation}

\section{Ergodic properties of simulated stochastic vector processes}
\label{A:2}

In this appendix, we prove that the sample functions generated by the third-order spectral representation method for stochastic vector processes are ergodic up to third-order.

Recall the simulation formula
\begin{equation}
\begin{aligned}
    f_{a}(t) &= 2 \sum_{k=0}^{N-1} \big[\sum_{m=1}^{m}|H_{am}(\omega_{mk})|\sqrt{\Delta \omega}\cos(\omega_{mk}t - \theta_{am}(\omega_{mk}) + \phi_{mk}) \\
    & + 2 \sum_{m=1}^{m}\sum_{n=1}^{m}\sum_{p=1}^{m}\sum_{q=1}^{m}\sum_{i + j = k}^{i \geq j\geq 0} |B_{amn}(\omega_{pi}, \omega_{qj})||G_{mp}(\omega_{pi})||G_{nq}(\omega_{qj})| \Delta \omega\\
    & \cos((\omega_{pi} + \omega_{qj})t - \beta_{amn}(\omega_{pi}, \omega_{qj}) - \theta^{I}_{mp}(\omega_{pi}) - \theta^{I}_{nq}(\omega_{qj}) + \phi_{pi} + \phi_{qj})\big]
\end{aligned}
\end{equation}
where the multi-indexed frequency is given by
\begin{equation}
\begin{aligned}
	& \omega_{mk} = k\Delta \omega + \frac{\Delta\omega}{2}(\frac{1}{N}) + \frac{\Delta\omega}{2}(\frac{m}{m})\\
\end{aligned}
\end{equation}
The smallest frequency wave corresponds to the largest period
\begin{equation}
\begin{aligned}
	& \omega_{01} = \Delta \omega + \frac{\Delta\omega}{2}(\frac{1}{N}) + \frac{\Delta\omega}{2}(\frac{1}{m}) = \frac{\Delta\omega}{2}(\frac{1}{m} + \frac{1}{N})\\
\end{aligned}
\end{equation}
yielding the wave period 
\begin{equation}
\begin{aligned}
	&T_{0} = \frac{4\pi mN}{m + N}\\
\end{aligned}
\end{equation}

\subsection{First Order Properties}
Ergodicity in the mean-value requires that:
\begin{equation}
    \langle f_a(t) \rangle_{T} = \mathbb{E}[f_a(t)] = 0
\end{equation}
Using the simulation equation above, we have
\begin{equation}
\begin{aligned}
	&\langle f_{a}(t) \rangle = \frac{2}{T}\int_{0}^{T}\sum_{k=0}^{N-1} \big[\sum_{m=1}^{m}|H_{am}(\omega_{mk})|\sqrt{\Delta \omega}\cos(\omega_{mk}t - \theta_{am}(\omega_{mk}) + \phi_{mk}) \\
    & + \sum_{m=1}^{m}\sum_{n=1}^{m}\sum_{p=1}^{m}\sum_{q=1}^{m}\sum_{i + j = k}^{i \geq j\geq 0} |B_{amn}(\omega_{pi}, \omega_{qj})||G_{mp}(\omega_{pi})||G_{nq}(\omega_{qj})| \Delta \omega\\
    & \cos((\omega_{pi} + \omega_{qj})t - \beta_{amn}(\omega_{pi}, \omega_{qj}) - \theta^{I}_{mp}(\omega_{pi}) - \theta^{I}_{nq}(\omega_{qj}) + \phi_{pi} + \phi_{qj})\big]dt\\
    & = 0\\
\end{aligned}
\end{equation}

\subsection{Second Order Properties}
Ergodicity in correlation requires that:
\begin{equation}
    \langle f_a(t) f_b(t + \tau) \rangle_{T} = E[f_a(t)f_b(t+\tau)] = R_{ab}(\tau)
\end{equation}
Using the proposed simulation formula, this can be expressed as
\begin{equation}
\begin{aligned}
	&\langle f_{a}(t)f_{b}(t + \tau) \rangle = \frac{4}{T}\int_{0}^{T} \sum_{k_{1}=0}^{N-1} \big[\sum_{l_{1}=1}^{m}|H_{al_{1}}(\omega_{l_{1}k_{1}})|\sqrt{\Delta \omega}\cos(\omega_{l_{1}k_{1}}t - \theta_{al_{1}}(\omega_{l_{1}k_{1}}) + \phi_{l_{1}k_{1}}) \\
    & + \sum_{l_{1}=1}^{m}\sum_{n_{1}=1}^{m}\sum_{p_{1}=1}^{m}\sum_{q_{1}=1}^{m}\sum_{i_{1} + j_{1} = k_{1}}^{i_{1} \geq j_{1} \geq 0} |B_{al_{1}n_{1}}(\omega_{p_{1}i_{1}}, \omega_{q_{1}j_{1}})||G_{l_{1}p_{1}}(\omega_{p_{1}i_{1}})||G_{n_{1}q_{1}}(\omega_{q_{1}j_{1}})| \Delta \omega\\
    & \cos((\omega_{p_{1}i_{1}} + \omega_{q_{1}j_{1}})t - \beta_{al_{1}n_{1}}(\omega_{p_{1}i_{1}}, \omega_{q_{1}j_{1}}) - \theta^{I}_{l_{1}p_{1}}(\omega_{p_{1}i_{1}}) - \theta^{I}_{n_{1}q_{1}}(\omega_{q_{1}j_{1}}) + \phi_{p_{1}i_{1}} + \phi_{q_{1}j_{1}})\big]\\
    & \sum_{k_{2}=0}^{N-1} \big[\sum_{l_{2}=1}^{m}|H_{al_{2}}(\omega_{l_{2}k_{2}})|\sqrt{\Delta \omega}\cos(\omega_{l_{2}k_{2}}(t + \tau) - \theta_{al_{2}}(\omega_{l_{2}k_{2}}) + \phi_{l_{2}k_{2}}) \\
    & + \sum_{l_{2}=1}^{m}\sum_{n_{2}=1}^{m}\sum_{p_{2}=1}^{m}\sum_{q_{2}=1}^{m}\sum_{i_{2} + j_{2} = k_{2}}^{i_{2} \geq j_{2} \geq 0} |B_{bl_{2}n_{2}}(\omega_{p_{2}i_{2}}, \omega_{q_{2}j_{2}})||G_{l_{2}p_{2}}(\omega_{p_{2}i_{2}})||G_{n_{2}q_{2}}(\omega_{q_{2}j_{2}})| \Delta \omega\\
    & \cos((\omega_{p_{2}i_{2}} + \omega_{q_{2}j_{2}})(t + \tau) - \beta_{bl_{2}n_{2}}(\omega_{p_{2}i_{2}}, \omega_{q_{2}j_{2}}) - \theta^{I}_{l_{2}p_{2}}(\omega_{p_{2}i_{2}}) - \theta^{I}_{n_{2}q_{2}}(\omega_{q_{2}j_{2}}) + \phi_{p_{2}i_{2}} + \phi_{q_{2}j_{2}})\big] dt\\
\end{aligned}
\end{equation}

\noindent
This product yields 4 terms. The first of these terms is given by
\begin{equation}
\begin{aligned}
	&\frac{4}{T}\int_{0}^{T} \sum_{k_{1}=0}^{N-1} \sum_{k_{2}=0}^{N-1} \sum_{l_{1}=1}^{m} \sum_{l_{2}=1}^{m} |H_{al_{1}}(\omega_{l_{1}k_{1}})||H_{bl_{2}}(\omega_{l_{2}k_{2}})| \Delta \omega \\
	& \cos(\omega_{l_{1}k_{1}}t - \theta_{al_{1}}(\omega_{l_{1}k_{1}}) + \phi_{l_{1}k_{1}}) \cos(\omega_{l_{2}k_{2}}(t + \tau) - \theta_{al_{2}}(\omega_{l_{2}k_{2}}) + \phi_{l_{2}k_{2}})\\
	& = \frac{2}{T}\int_{0}^{T} \sum_{k_{1}=0}^{N-1} \sum_{k_{2}=0}^{N-1} \sum_{l_{1}=1}^{m} \sum_{l_{2}=1}^{m} |H_{al_{1}}(\omega_{l_{1}k_{1}})||H_{bl_{2}}(\omega_{l_{2}k_{2}})| \Delta \omega \\
	& \Big[ \cos((\omega_{l_{1}k_{1}} + \omega_{l_{2}k_{2}})t + \omega_{l_{2}k_{2}}\tau - \theta_{al_{1}}(\omega_{l_{1}k_{1}}) - \theta_{bl_{2}}(\omega_{l_{2}k_{2}}) + \phi_{l_{1}k_{1}} + \phi_{l_{2}k_{2}}) \\
	& \cos((\omega_{l_{2}k_{2}} - \omega_{l_{1}k_{1}})t + \omega_{l_{2}k_{2}}\tau + \theta_{al_{1}}(\omega_{l_{1}k_{1}}) - \theta_{bl_{2}}(\omega_{l_{2}k_{2}}) - \phi_{l_{1}k_{1}} + \phi_{l_{2}k_{2}}) \Big]dt \\
\end{aligned}
\end{equation}
Recognizing that
\begin{equation}
\begin{aligned}
	&\int_{0}^{T} \cos((\omega_{l_{1}k_{1}} + \omega_{l_{2}k_{2}})t + \omega_{l_{2}k_{2}}\tau - \theta_{al_{1}}(\omega_{l_{1}k_{1}}) - \theta_{al_{2}}(\omega_{l_{2}k_{2}}) + \phi_{l_{1}k_{1}} + \phi_{l_{2}k_{2}}) dt = 0; \\ & \forall l_{1}, l_{2}, k_{1}, k_{2}\\
	&\int_{0}^{T} \cos((\omega_{l_{2}k_{2}} - \omega_{l_{1}k_{1}})t + \omega_{l_{2}k_{2}}\tau + \theta_{al_{1}}(\omega_{l_{1}k_{1}}) - \theta_{al_{2}}(\omega_{l_{2}k_{2}}) - \phi_{l_{1}k_{1}} + \phi_{l_{2}k_{2}}) dt = 0\\
	& \forall l_{1}, l_{2}, k_{1}, k_{2}; \ l_{1} \neq l_{2}, k_{1} \neq k_{2}\\
\end{aligned}
\end{equation}
we have
\begin{equation}
\begin{aligned}
	\langle f_{a}(t)f_{b}(t + \tau) \rangle & = \frac{2}{T} \int_{0}^{T} \sum_{k=0}^{N-1}  \sum_{m=1}^{m} |H_{am}(\omega_{mk})||H_{bm}(\omega_{mk})| \Delta \omega \cos(\omega_{mk}\tau + \theta_{am}(\omega_{mk}) - \theta_{bm}(\omega_{mk})) dt \\
	& = \frac{2}{T} \int_{0}^{T} \sum_{k=0}^{N-1} S_{ab}^{P}(\omega_{k}) e^{\iota \omega_{k} \tau} \Delta \omega dt \\
	& = 2 \sum_{k=0}^{N-1} S_{ab}^{P}(\omega_{k}) e^{\iota \omega_{k} \tau} \Delta \omega \\
\end{aligned}
\end{equation}

\noindent
The second and the third term are given by
\begin{equation}
\begin{aligned}
	&\frac{4}{T}\int_{0}^{T} \sum_{k_{1}=0}^{N-1} \sum_{k_{2}=0}^{N-1} \sum_{l_{1}=1}^{m} \sum_{l_{2}=1}^{m} \sum_{n_{2}=1}^{m} \sum_{p_{2}=1}^{m} \sum_{q_{2}=1}^{m} \sum_{i_{2} + j_{2} = k_{2}}^{i_{2} \geq j_{2} \geq 0} \\
	& |H_{al_{1}}(\omega_{l_{1}k_{1}})| |B_{bl_{2}n_{2}}(\omega_{p_{2}i_{2}}, \omega_{q_{2}j_{2}})||G_{l_{2}p_{2}}(\omega_{p_{2}i_{2}})||G_{n_{2}q_{2}}(\omega_{q_{2}j_{2}})| {\Delta \omega}^{\frac{3}{2}} \\
	& \cos(\omega_{l_{1}k_{1}}t - \theta_{al_{1}}(\omega_{l_{1}k_{1}}) + \phi_{l_{1}k_{1}})\\
    & \cos((\omega_{p_{2}i_{2}} + \omega_{q_{2}j_{2}})(t + \tau) - \beta_{bl_{2}n_{2}}(\omega_{p_{2}i_{2}}, \omega_{q_{2}j_{2}}) - \theta^{I}_{l_{2}p_{2}}(\omega_{p_{2}i_{2}}) - \theta^{I}_{n_{2}q_{2}}(\omega_{q_{2}j_{2}}) + \phi_{p_{2}i_{2}} + \phi_{q_{2}j_{2}})dt \\
    & = \frac{2}{T}\int_{0}^{T} \sum_{k_{1}=0}^{N-1} \sum_{k_{2}=0}^{N-1} \sum_{l_{1}=1}^{m} \sum_{l_{2}=1}^{m} \sum_{n_{2}=1}^{m} \sum_{p_{2}=1}^{m} \sum_{q_{2}=1}^{m} \sum_{i_{2} + j_{2} = k_{2}}^{i_{2} \geq j_{2} \geq 0} \\
	& |H_{al_{1}}(\omega_{l_{1}k_{1}})| |B_{bl_{2}n_{2}}(\omega_{p_{2}i_{2}}, \omega_{q_{2}j_{2}})||G_{l_{2}p_{2}}(\omega_{p_{2}i_{2}})||G_{n_{2}q_{2}}(\omega_{q_{2}j_{2}})| {\Delta \omega}^{\frac{3}{2}} \\
	& \Big[ \cos((\omega_{l_{1}k_{1}} + \omega_{p_{2}i_{2}} + \omega_{q_{2}j_{2}})t + (\omega_{p_{2}i_{2}} + \omega_{q_{2}j_{2}})\tau \\
	& - \theta_{al_{1}}(\omega_{l_{1}k_{1}}) - \beta_{bl_{2}n_{2}}(\omega_{p_{2}i_{2}}, \omega_{q_{2}j_{2}}) - \theta^{I}_{l_{2}p_{2}}(\omega_{p_{2}i_{2}}) - \theta^{I}_{n_{2}q_{2}}(\omega_{q_{2}j_{2}}) + \phi_{l_{1}k_{1}} + \phi_{p_{2}i_{2}} + \phi_{q_{2}j_{2}}) \\
	& \cos((\omega_{p_{2}i_{2}} + \omega_{q_{2}j_{2}} - \omega_{l_{1}k_{1}})t + (\omega_{p_{2}i_{2}} + \omega_{q_{2}j_{2}})\tau \\
	& + \theta_{al_{1}}(\omega_{l_{1}k_{1}}) - \beta_{bl_{2}n_{2}}(\omega_{p_{2}i_{2}}, \omega_{q_{2}j_{2}}) - \theta^{I}_{l_{2}p_{2}}(\omega_{p_{2}i_{2}}) - \theta^{I}_{n_{2}q_{2}}(\omega_{q_{2}j_{2}}) - \phi_{l_{1}k_{1}} + \phi_{p_{2}i_{2}} + \phi_{q_{2}j_{2}}) \Big]dt \\
	& = 0\\
\end{aligned}
\end{equation}
because
\begin{equation}
\begin{aligned}
	&\int_{0}^{T} \cos((\omega_{l_{1}k_{1}} + \omega_{p_{2}i_{2}} + \omega_{q_{2}j_{2}})t) dt = 0; \ \forall l_{1}, p_{2}, q_{2}, k_{1}, i_{2}, j_{2} \\
	&\int_{0}^{T} \cos((\omega_{p_{2}i_{2}} + \omega_{q_{2}j_{2}} - \omega_{l_{1}k_{1}})t) dt = 0; \ \forall l_{1}, p_{2}, q_{2}, k_{1}, i_{2}, j_{2} \\
\end{aligned}
\end{equation}

Finally, the fourth term is given by
\begin{equation}
\begin{aligned}
    & \frac{4}{T}\int_{0}^{T} \sum_{k_{1}=0}^{N-1} \sum_{k_{2}=0}^{N-1} \sum_{l_{1}=1}^{m} \sum_{n_{1}=1}^{m} \sum_{p_{1}=1}^{m} \sum_{q_{1}=1}^{m} \sum_{i_{1} + j_{1} = k_{1}}^{i_{1} \geq j_{1} \geq 0} \sum_{l_{2}=1}^{m} \sum_{n_{2}=1}^{m} \sum_{p_{2}=1}^{m} \sum_{q_{2}=1}^{m} \sum_{i_{2} + j_{2} = k_{2}}^{i_{2} \geq j_{2} \geq 0} \\
    & |B_{al_{1}n_{1}}(\omega_{p_{1}i_{1}}, \omega_{q_{1}j_{1}})||G_{l_{1}p_{1}}(\omega_{p_{1}i_{1}})||G_{n_{1}q_{1}}(\omega_{q_{1}j_{1}})| \\
    & |B_{bl_{2}n_{2}}(\omega_{p_{2}i_{2}}, \omega_{q_{2}j_{2}})||G_{l_{2}p_{2}}(\omega_{p_{2}i_{2}})||G_{n_{2}q_{2}}(\omega_{q_{2}j_{2}})| {\Delta \omega}^{2} \\
    & \cos((\omega_{p_{1}i_{1}} + \omega_{q_{1}j_{1}})t - \beta_{al_{1}n_{1}}(\omega_{p_{1}i_{1}}, \omega_{q_{1}j_{1}}) - \theta^{I}_{l_{1}p_{1}}(\omega_{p_{1}i_{1}}) - \theta^{I}_{n_{1}q_{1}}(\omega_{q_{1}j_{1}}) + \phi_{p_{1}i_{1}} + \phi_{q_{1}j_{1}}) \\
    & \cos((\omega_{p_{2}i_{2}} + \omega_{q_{2}j_{2}})(t + \tau) - \beta_{bl_{2}n_{2}}(\omega_{p_{2}i_{2}}, \omega_{q_{2}j_{2}}) - \theta^{I}_{l_{2}p_{2}}(\omega_{p_{2}i_{2}}) - \theta^{I}_{n_{2}q_{2}}(\omega_{q_{2}j_{2}}) + \phi_{p_{2}i_{2}} + \phi_{q_{2}j_{2}})dt \\
    & = \frac{4}{T}\int_{0}^{T} \sum_{k_{1}=0}^{N-1} \sum_{k_{2}=0}^{N-1} \sum_{l_{1}=1}^{m} \sum_{n_{1}=1}^{m} \sum_{p_{1}=1}^{m} \sum_{q_{1}=1}^{m} \sum_{i_{1} + j_{1} = k_{1}}^{i_{1} \geq j_{1} \geq 0} \sum_{l_{2}=1}^{m} \sum_{n_{2}=1}^{m} \sum_{p_{2}=1}^{m} \sum_{q_{2}=1}^{m} \sum_{i_{2} + j_{2} = k_{2}}^{i_{2} \geq j_{2} \geq 0} \\
    & |B_{al_{1}n_{1}}(\omega_{p_{1}i_{1}}, \omega_{q_{1}j_{1}})||G_{l_{1}p_{1}}(\omega_{p_{1}i_{1}})||G_{n_{1}q_{1}}(\omega_{q_{1}j_{1}})| \\
    & |B_{bl_{2}n_{2}}(\omega_{p_{2}i_{2}}, \omega_{q_{2}j_{2}})||G_{l_{2}p_{2}}(\omega_{p_{2}i_{2}})||G_{n_{2}q_{2}}(\omega_{q_{2}j_{2}})| {\Delta \omega}^{2} \\
    & \cos((\omega_{p_{1}i_{1}} + \omega_{q_{1}j_{1}})t - \beta_{al_{1}n_{1}}(\omega_{p_{1}i_{1}}, \omega_{q_{1}j_{1}}) - \theta^{I}_{l_{1}p_{1}}(\omega_{p_{1}i_{1}}) - \theta^{I}_{n_{1}q_{1}}(\omega_{q_{1}j_{1}}) + \phi_{p_{1}i_{1}} + \phi_{q_{1}j_{1}}) \\
    & \cos((\omega_{p_{2}i_{2}} + \omega_{q_{2}j_{2}})(t + \tau) - \beta_{bl_{2}n_{2}}(\omega_{p_{2}i_{2}}, \omega_{q_{2}j_{2}}) - \theta^{I}_{l_{2}p_{2}}(\omega_{p_{2}i_{2}}) - \theta^{I}_{n_{2}q_{2}}(\omega_{q_{2}j_{2}}) + \phi_{p_{2}i_{2}} + \phi_{q_{2}j_{2}})dt \\
    & = \frac{2}{T}\int_{0}^{T} \sum_{k_{1}=0}^{N-1} \sum_{k_{2}=0}^{N-1} \sum_{l_{1}=1}^{m} \sum_{n_{1}=1}^{m} \sum_{p_{1}=1}^{m} \sum_{q_{1}=1}^{m} \sum_{i_{1} + j_{1} = k_{1}}^{i_{1} \geq j_{1} \geq 0} \sum_{l_{2}=1}^{m} \sum_{n_{2}=1}^{m} \sum_{p_{2}=1}^{m} \sum_{q_{2}=1}^{m} \sum_{i_{2} + j_{2} = k_{2}}^{i_{2} \geq j_{2} \geq 0} \\
    & |B_{al_{1}n_{1}}(\omega_{p_{1}i_{1}}, \omega_{q_{1}j_{1}})||G_{l_{1}p_{1}}(\omega_{p_{1}i_{1}})||G_{n_{1}q_{1}}(\omega_{q_{1}j_{1}})| \\
    & |B_{bl_{2}n_{2}}(\omega_{p_{2}i_{2}}, \omega_{q_{2}j_{2}})||G_{l_{2}p_{2}}(\omega_{p_{2}i_{2}})||G_{n_{2}q_{2}}(\omega_{q_{2}j_{2}})| {\Delta \omega}^{2} \\
    & \Big[ \cos((\omega_{p_{1}i_{1}} + \omega_{q_{1}j_{1}} + \omega_{p_{2}i_{2}} + \omega_{q_{2}j_{2}})t + (\omega_{p_{2}i_{2}} + \omega_{q_{2}j_{2}})\tau - \beta_{al_{1}n_{1}}(\omega_{p_{1}i_{1}}, \omega_{q_{1}j_{1}}) - \beta_{bl_{2}n_{2}}(\omega_{p_{2}i_{2}}, \omega_{q_{2}j_{2}}) \\
    & - \theta^{I}_{l_{1}p_{1}}(\omega_{p_{1}i_{1}}) - \theta^{I}_{n_{1}q_{1}}(\omega_{q_{1}j_{1}}) - \theta^{I}_{l_{2}p_{2}}(\omega_{p_{2}i_{2}}) - \theta^{I}_{n_{2}q_{2}}(\omega_{q_{2}j_{2}}) + \phi_{p_{1}i_{1}} + \phi_{q_{1}j_{1}} + \phi_{p_{2}i_{2}} + \phi_{q_{2}j_{2}}) \\
    & \cos((\omega_{p_{2}i_{2}} + \omega_{q_{2}j_{2}} - \omega_{p_{1}i_{1}} - \omega_{q_{1}j_{1}})t + (\omega_{p_{2}i_{2}} + \omega_{q_{2}j_{2}})\tau + \beta_{al_{1}n_{1}}(\omega_{p_{1}i_{1}}, \omega_{q_{1}j_{1}}) - \beta_{bl_{2}n_{2}}(\omega_{p_{2}i_{2}}, \omega_{q_{2}j_{2}}) \\
    & + \theta^{I}_{l_{1}p_{1}}(\omega_{p_{1}i_{1}}) + \theta^{I}_{n_{1}q_{1}}(\omega_{q_{1}j_{1}}) - \theta^{I}_{l_{2}p_{2}}(\omega_{p_{2}i_{2}}) - \theta^{I}_{n_{2}q_{2}}(\omega_{q_{2}j_{2}}) - \phi_{p_{1}i_{1}} - \phi_{q_{1}j_{1}} + \phi_{p_{2}i_{2}} + \phi_{q_{2}j_{2}}) \Big]dt \\
\end{aligned}
\end{equation}
Recognizing that
\begin{equation}
\begin{aligned}
	&\int_{0}^{T} \cos((\omega_{p_{1}i_{1}} + \omega_{q_{1}j_{1}} + \omega_{p_{2}i_{2}} + \omega_{q_{2}j_{2}})t) dt = 0; \ \forall p_{1}, q_{1}, p_{2}, q_{2}, i_{1}, j_{1}, i_{2}, j_{2} \\
	&\int_{0}^{T} \cos((\omega_{p_{2}i_{2}} + \omega_{q_{2}j_{2}} - \omega_{p_{1}i_{1}} - \omega_{q_{1}j_{1}})t) dt = 0; \ \forall p_{1}, q_{1}, p_{2}, q_{2}, i_{1}, j_{1}, i_{2}, j_{2}; \ p_{1} \neq p_{2}, q_{1} \neq q_{2}, i_{1} \neq i_{2}, j_{1} \neq j_{2} \\
\end{aligned}
\end{equation}
this expression simplifies as
\begin{equation}
\begin{aligned}
	& = \frac{2}{T}\int_{0}^{T} \sum_{k=0}^{N-1} \sum_{l_{1}=1}^{m} \sum_{n_{1}=1}^{m} \sum_{l_{2}=1}^{m} \sum_{n_{2}=1}^{m} \sum_{p=1}^{m} \sum_{q=1}^{m} \sum_{i + j = k}^{i \geq j \geq 0} |B_{al_{1}n_{1}}(\omega_{pi}, \omega_{qj})||G_{l_{1}p}(\omega_{pi})||G_{n_{1}q}(\omega_{qj})| \\
    & |B_{bl_{2}n_{2}}(\omega_{pi}, \omega_{qj})||G_{l_{2}p}(\omega_{pi})||G_{n_{2}q}(\omega_{qj})| {\Delta \omega}^{2} \\
    & \cos((\omega_{pi} + \omega_{qj})\tau + \beta_{al_{1}n_{1}}(\omega_{pi}, \omega_{qj}) - \beta_{bl_{2}n_{2}}(\omega_{pi}, \omega_{qj}) + \theta^{I}_{l_{1}p}(\omega_{pi}) + \theta^{I}_{n_{1}q}(\omega_{qj}) - \theta^{I}_{l_{2}p}(\omega_{pi}) - \theta^{I}_{n_{2}q}(\omega_{qj}))dt \\
    & = \frac{2}{T}\int_{0}^{T} \sum_{k=0}^{N-1} \sum_{l_{1}=1}^{m} \sum_{n_{1}=1}^{m} \sum_{l_{2}=1}^{m} \sum_{n_{2}=1}^{m} \sum_{p=1}^{m} \sum_{q=1}^{m} \sum_{i + j = k}^{i \geq j \geq 0} B_{al_{1}n_{1}}(\omega_{pi}, \omega_{qj})G_{l_{1}p}(\omega_{pi})G_{n_{1}q}(\omega_{qj})\\
    & B^{*}_{bl_{2}n_{2}}(\omega_{pi}, \omega_{qj})G^{*}_{l_{2}p}(\omega_{pi})G^{*}_{n_{2}q}(\omega_{qj}) {\Delta \omega}^{2} \cos((\omega_{pi} + \omega_{qj})\tau)dt \\
    & = \frac{2}{T}\int_{0}^{T} \sum_{k=0}^{N-1} S^{I}_{ab}(\omega_{k}) \Delta \omega \cos((\omega_{i} + \omega_{j})\tau)dt \\
    & = \frac{2}{T}\int_{0}^{T} \sum_{k=0}^{N-1} S^{I}_{ab}(\omega_{k}) \Delta \omega \cos(\omega_{k}\tau)dt \\
    & = 2 \sum_{k=0}^{N-1} S^{I}_{ab}(\omega_{k}) \Delta \omega \cos(\omega_{k}\tau) \\
\end{aligned}
\end{equation}

\noindent
Combining these terms leads to
\begin{equation}
\begin{aligned}
    \langle f_{a}(t)f_{b}(t + \tau) \rangle &= 2 \sum_{k=0}^{N-1} S_{ab}(\omega_{k})\Delta \omega \cos(\omega_{k}\tau) \\
    & = R_{ab}(\tau)\\
\end{aligned}
\end{equation}

\subsection{Third Order Properties}

Ergodicity in the third-order moment requires that:
\begin{equation}
    \langle f_a(t) f_b(t + \tau_1) f_c(t+\tau_2) \rangle_{T} = E[f_a(t)f_b(t+\tau_1)f_c(t+\tau_2)] = R^{(3)}(\tau_1,\tau_2)
\end{equation}
Using the proposed simulation formula, this can be expressed as
\begin{equation}
\begin{aligned}
	&\langle f_{a}(t)f_{b}(t + \tau_{1})f_{c}(t + \tau_{2}) \rangle \\
	& = \frac{8}{T}\int_{0}^{T} \sum_{k_{1}=0}^{N-1} \big[\sum_{l_{1}=1}^{m}|H_{al_{1}}(\omega_{l_{1}k_{1}})|\sqrt{\Delta \omega}\cos(\omega_{l_{1}k_{1}}t - \theta_{al_{1}}(\omega_{l_{1}k_{1}}) + \phi_{l_{1}k_{1}}) \\
    & + \sum_{l_{1}=1}^{m}\sum_{n_{1}=1}^{m}\sum_{p_{1}=1}^{m}\sum_{q_{1}=1}^{m}\sum_{i_{1} + j_{1} = k_{1}}^{i_{1} \geq j_{1} \geq 0} |B_{al_{1}n_{1}}(\omega_{p_{1}i_{1}}, \omega_{q_{1}j_{1}})||G_{l_{1}p_{1}}(\omega_{p_{1}i_{1}})||G_{n_{1}q_{1}}(\omega_{q_{1}j_{1}})| \Delta \omega\\
    & \cos((\omega_{p_{1}i_{1}} + \omega_{q_{1}j_{1}})t - \beta_{al_{1}n_{1}}(\omega_{p_{1}i_{1}}, \omega_{q_{1}j_{1}}) - \theta^{I}_{l_{1}p_{1}}(\omega_{p_{1}i_{1}}) - \theta^{I}_{n_{1}q_{1}}(\omega_{q_{1}j_{1}}) + \phi_{p_{1}i_{1}} + \phi_{q_{1}j_{1}})\big]\\
    & \sum_{k_{2}=0}^{N-1} \big[\sum_{l_{2}=1}^{m}|H_{al_{2}}(\omega_{l_{2}k_{2}})|\sqrt{\Delta \omega}\cos(\omega_{l_{2}k_{2}}(t + \tau_{1}) - \theta_{al_{2}}(\omega_{l_{2}k_{2}}) + \phi_{l_{2}k_{2}}) \\
    & + \sum_{l_{2}=1}^{m}\sum_{n_{2}=1}^{m}\sum_{p_{2}=1}^{m}\sum_{q_{2}=1}^{m}\sum_{i_{2} + j_{2} = k_{2}}^{i_{2} \geq j_{2} \geq 0} |B_{bl_{2}n_{2}}(\omega_{p_{2}i_{2}}, \omega_{q_{2}j_{2}})||G_{l_{2}p_{2}}(\omega_{p_{2}i_{2}})||G_{n_{2}q_{2}}(\omega_{q_{2}j_{2}})| \Delta \omega\\
    & \cos((\omega_{p_{2}i_{2}} + \omega_{q_{2}j_{2}})(t + \tau_{1}) - \beta_{bl_{2}n_{2}}(\omega_{p_{2}i_{2}}, \omega_{q_{2}j_{2}}) - \theta^{I}_{l_{2}p_{2}}(\omega_{p_{2}i_{2}}) - \theta^{I}_{n_{2}q_{2}}(\omega_{q_{2}j_{2}}) + \phi_{p_{2}i_{2}} + \phi_{q_{2}j_{2}})\big]\\
    & \sum_{k_{3}=0}^{N-1} \big[\sum_{l_{3}=1}^{m}|H_{al_{3}}(\omega_{l_{3}k_{3}})|\sqrt{\Delta \omega}\cos(\omega_{l_{3}k_{3}}(t + \tau_{2}) - \theta_{al_{3}}(\omega_{l_{3}k_{3}}) + \phi_{l_{3}k_{3}}) \\
    & + \sum_{l_{3}=1}^{m}\sum_{n_{3}=1}^{m}\sum_{p_{3}=1}^{m}\sum_{q_{3}=1}^{m}\sum_{i_{3} + j_{3} = k_{3}}^{i_{3} \geq j_{3} \geq 0} |B_{bl_{3}n_{3}}(\omega_{p_{3}i_{3}}, \omega_{q_{3}j_{3}})||G_{l_{3}p_{3}}(\omega_{p_{3}i_{3}})||G_{n_{3}q_{3}}(\omega_{q_{3}j_{3}})| \Delta \omega\\
    & \cos((\omega_{p_{3}i_{3}} + \omega_{q_{3}j_{3}})(t + \tau_{2}) - \beta_{bl_{3}n_{3}}(\omega_{p_{3}i_{3}}, \omega_{q_{3}j_{3}}) - \theta^{I}_{l_{3}p_{3}}(\omega_{p_{3}i_{3}}) - \theta^{I}_{n_{3}q_{3}}(\omega_{q_{3}j_{3}}) + \phi_{p_{3}i_{3}} + \phi_{q_{3}j_{3}})\big] dt\\
\end{aligned}
\end{equation}

\begin{equation}
\begin{aligned}
    & \langle f_{a}(t)f_{b}(t + \tau_{1})f_{b}(t + \tau_{2}) \rangle \\
    & = \frac{8}{T}\int_{0}^{T} \sum_{k_{2}=0}^{N-1} \sum_{l_{2}=1}^{m}|H_{bl_{2}}(\omega_{l_{2}k_{2}})|\sqrt{\Delta \omega}\cos(\omega_{l_{2}k_{2}}(t + \tau_{1}) - \theta_{bl_{2}}(\omega_{l_{2}k_{2}}) + \phi_{l_{2}k_{2}}) \\
    & \sum_{k_{3}=0}^{N-1} \sum_{l_{3}=1}^{m}|H_{cl_{3}}(\omega_{l_{3}k_{3}})|\sqrt{\Delta \omega}\cos(\omega_{l_{3}k_{3}}(t + \tau_{2}) - \theta_{cl_{3}}(\omega_{l_{3}k_{3}}) + \phi_{l_{3}k_{3}})\\
    & \sum_{k_{1}=0}^{N-1}\sum_{l_{1}=1}^{m}\sum_{n_{1}=1}^{m}\sum_{p_{1}=1}^{m}\sum_{q_{1}=1}^{m}\sum_{i_{1} + j_{1} = k_{1}}^{i_{1} \geq j_{1} \geq 0} |B_{al_{1}n_{1}}(\omega_{p_{1}i_{1}}, \omega_{q_{1}j_{1}})||G_{l_{1}p_{1}}(\omega_{p_{1}i_{1}})||G_{n_{1}q_{1}}(\omega_{q_{1}j_{1}})| \Delta \omega \\
    & \cos((\omega_{p_{1}i_{1}} + \omega_{q_{1}j_{1}})t - \beta_{al_{1}n_{1}}(\omega_{p_{1}i_{1}}, \omega_{q_{1}j_{1}}) - \theta^{I}_{l_{1}p_{1}}(\omega_{p_{1}i_{1}}) - \theta^{I}_{n_{1}q_{1}}(\omega_{q_{1}j_{1}}) + \phi_{p_{1}i_{1}} + \phi_{q_{1}j_{1}}) dt \\
    & = \frac{24}{T}\int_{0}^{T} \sum_{k_{1}=0}^{N-1}\sum_{l_{1}=1}^{m}\sum_{n_{1}=1}^{m}\sum_{p_{1}=1}^{m}\sum_{q_{1}=1}^{m}\sum_{i_{1} + j_{1} = k_{1}}^{i_{1} \geq j_{1} \geq 0}\sum_{k_{2}=0}^{N-1} \sum_{l_{2}=1}^{m}\sum_{k_{3}=0}^{N-1}\sum_{l_{3}=1}^{m} \\
    & |B_{al_{1}n_{1}}(\omega_{p_{1}i_{1}}, \omega_{q_{1}j_{1}})||G_{l_{1}p_{1}}(\omega_{p_{1}i_{1}})||G_{n_{1}q_{1}}(\omega_{q_{1}j_{1}})||H_{bl_{2}}(\omega_{l_{2}k_{2}})||H_{cl_{3}}(\omega_{l_{3}k_{3}})| (\Delta \omega)^{2} \\
    & \cos(\omega_{l_{2}k_{2}}(t +\tau_{1} ) - \theta_{bl_{2}}(\omega_{l_{2}k_{2}}) + \phi_{l_{2}k_{2}})\cos(\omega_{l_{3}k_{3}}(t + \tau_{2}) - \theta_{cl_{3}}(\omega_{l_{3}k_{3}}) + \phi_{l_{3}k_{3}}) \\
    & \cos((\omega_{p_{1}i_{1}} + \omega_{q_{1}j_{1}})t - \beta_{al_{1}n_{1}}(\omega_{p_{1}i_{1}}, \omega_{q_{1}j_{1}}) - \theta^{I}_{l_{1}p_{1}}(\omega_{p_{1}i_{1}}) - \theta^{I}_{n_{1}q_{1}}(\omega_{q_{1}j_{1}}) + \phi_{p_{1}i_{1}} + \phi_{q_{1}j_{1}}) dt \\
    & = \frac{12}{T}\int_{0}^{T} \sum_{k_{1}=0}^{N-1}\sum_{l_{1}=1}^{m}\sum_{n_{1}=1}^{m}\sum_{p_{1}=1}^{m}\sum_{q_{1}=1}^{m}\sum_{i_{1} + j_{1} = k_{1}}^{i_{1} \geq j_{1} \geq 0}\sum_{k_{2}=0}^{N-1} \sum_{l_{2}=1}^{m}\sum_{k_{3}=0}^{N-1}\sum_{l_{3}=1}^{m} \\
    & |B_{al_{1}n_{1}}(\omega_{p_{1}i_{1}}, \omega_{q_{1}j_{1}})||G_{l_{1}p_{1}}(\omega_{p_{1}i_{1}})||G_{n_{1}q_{1}}(\omega_{q_{1}j_{1}})||H_{bl_{2}}(\omega_{l_{2}k_{2}})||H_{cl_{3}}(\omega_{l_{3}k_{3}})| \Delta \omega \\
    & \Big[ \cos((\omega_{l_{2}k_{2}} + \omega_{l_{3}k_{3}})t + \omega_{l_{2}k_{2}}\tau_{1} + \omega_{l_{3}k_{3}}\tau_2 - \theta_{bl_{2}}(\omega_{l_{2}k_{2}}) - \theta_{cl_{3}}(\omega_{l_{3}k_{3}}) + \phi_{l_{2}k_{2}} + \phi_{l_{3}k_{3}}) \\
    & + \cos((\omega_{l_{2}k_{2}} - \omega_{l_{3}k_{3}})t + \omega_{l_{2}k_{2}}\tau_{1} - \omega_{l_{3}k_{3}}\tau_2 - \theta_{bl_{2}}(\omega_{l_{2}k_{2}}) + \theta_{cl_{3}}(\omega_{l_{3}k_{3}}) + \phi_{l_{2}k_{2}} - \phi_{l_{3}k_{3}})\Big] \\
    & \cos((\omega_{p_{1}i_{1}} + \omega_{q_{1}j_{1}})t - \beta_{al_{1}n_{1}}(\omega_{p_{1}i_{1}}, \omega_{q_{1}j_{1}}) - \theta^{I}_{l_{1}p_{1}}(\omega_{p_{1}i_{1}}) - \theta^{I}_{n_{1}q_{1}}(\omega_{q_{1}j_{1}}) + \phi_{p_{1}i_{1}} + \phi_{q_{1}j_{1}}) dt \\
\end{aligned}
\end{equation}

the expression further simplifies to

\begin{equation}
\begin{aligned}
	& = \frac{12}{T}\int_{0}^{T} \sum_{k_{1}=0}^{N-1}\sum_{l_{1}=1}^{m}\sum_{n_{1}=1}^{m}\sum_{p_{1}=1}^{m}\sum_{q_{1}=1}^{m}\sum_{i_{1} + j_{1} = k_{1}}^{i_{1} \geq j_{1} \geq 0}\sum_{k_{2}=0}^{N-1} \sum_{l_{2}=1}^{m}\sum_{k_{3}=0}^{N-1}\sum_{l_{3}=1}^{m} \\
    & |B_{al_{1}n_{1}}(\omega_{p_{1}i_{1}}, \omega_{q_{1}j_{1}})||G_{l_{1}p_{1}}(\omega_{p_{1}i_{1}})||G_{n_{1}q_{1}}(\omega_{q_{1}j_{1}})||H_{bl_{2}}(\omega_{l_{2}k_{2}})||H_{cl_{3}}(\omega_{l_{3}k_{3}})| (\Delta \omega)^{2} \\
    & \Big[ \cos((\omega_{l_{2}k_{2}} + \omega_{l_{3}k_{3}})t + \omega_{l_{2}k_{2}}\tau_{1} + \omega_{l_{3}k_{3}}\tau_2 - \theta_{bl_{2}}(\omega_{l_{2}k_{2}}) - \theta_{cl_{3}}(\omega_{l_{3}k_{3}}) + \phi_{l_{2}k_{2}} + \phi_{l_{3}k_{3}})\\
    & \cos((\omega_{p_{1}i_{1}} + \omega_{q_{1}j_{1}})t - \beta_{al_{1}n_{1}}(\omega_{p_{1}i_{1}}, \omega_{q_{1}j_{1}}) - \theta^{I}_{l_{1}p_{1}}(\omega_{p_{1}i_{1}}) - \theta^{I}_{n_{1}q_{1}}(\omega_{q_{1}j_{1}}) + \phi_{p_{1}i_{1}} + \phi_{q_{1}j_{1}}) \\
    & + \cos((\omega_{l_{2}k_{2}} - \omega_{l_{3}k_{3}})t + \omega_{l_{2}k_{2}}\tau_{1} - \omega_{l_{3}k_{3}}\tau_2 - \theta_{bl_{2}}(\omega_{l_{2}k_{2}}) + \theta_{cl_{3}}(\omega_{l_{3}k_{3}}) + \phi_{l_{2}k_{2}} - \phi_{l_{3}k_{3}}) \\
    & \cos((\omega_{p_{1}i_{1}} + \omega_{q_{1}j_{1}})t - \beta_{al_{1}n_{1}}(\omega_{p_{1}i_{1}}, \omega_{q_{1}j_{1}}) - \theta^{I}_{l_{1}p_{1}}(\omega_{p_{1}i_{1}}) - \theta^{I}_{n_{1}q_{1}}(\omega_{q_{1}j_{1}}) + \phi_{p_{1}i_{1}} + \phi_{q_{1}j_{1}}) \Big] dt \\
\end{aligned}
\end{equation}

Of the two sets of cosine product terms above, the second product sums to zero because of odd number of random phase angles in the expansion. Simplifying the first cosine product more and recognizing that

\begin{equation}
\begin{aligned}
	&\int_{0}^{T} \cos(\omega_{l_{2}k_{2}}t + \omega_{l_{3}k_{3}}t - (\omega_{p_{1}i_{1}} + \omega_{q_{1}j_{1}})t) dt = 0; \ \forall p_{1}, q_{1}, p_{2}, q_{2}, i_{1}, j_{1}, i_{2}, j_{2}; \ p_{1} \neq l_{2}, q_{1} \neq l_{3}, i_{1} \neq k_{2}, j_{1} \neq k_{3} \\
\end{aligned}
\end{equation}

and

\begin{equation}
\begin{aligned}
	&G_{ab}(\omega)H_{bc}(\omega) = I_{ac}\\
\end{aligned}
\end{equation}

it leads to

\begin{equation}
\begin{aligned}
	& = \frac{6}{T}\int_{0}^{T} \sum_{k_{2}=0}^{N-1} \sum_{k_{3}=0}^{N-1} B_{abc}(\omega_{k_{2}}, \omega_{k_{3}}) (\Delta \omega)^{2}\cos(\omega_{k_{2}}\tau_{1} + \omega_{k_{3}}\tau_{2}) dt\\
	& = 6 \sum_{k_{2}=0}^{N-1} \sum_{k_{3}=0}^{N-1} B_{abc}(\omega_{k_{2}}, \omega_{k_{3}}) (\Delta \omega)^{2}\cos(\omega_{k_{2}}\tau_{1} + \omega_{k_{3}}\tau_{2}) \\
	& = R^{(3)}_{abc}(\tau_{1}, \tau_{2}) \\
\end{aligned}
\end{equation}

\end{document}